\documentclass[11pt]{amsart}
\usepackage{amsmath,amsfonts,amssymb,amsthm,epsfig,verbatim,ifthen,graphicx}
\usepackage[usenames,dvips]{color}
\usepackage[all]{xy}
\theoremstyle{definition} 
\newtheorem{definition}{Definition}[section]
\newtheorem{theorem}[definition]{Theorem}
\newtheorem{proposition}[definition]{Proposition}
\newtheorem{lemma}[definition]{Lemma}

\theoremstyle{definition}
\newtheorem{corollary}[definition]{Corollary}

\newtheorem{remarks}[definition]{Remarks}
\newtheorem{example}[definition]{Example}
\newtheorem{remark}[definition]{Remark}

\newcommand{\N}{\mathbb{N}}
\newcommand{\Z}{\mathbb{Z}}
\newcommand{\R}{\mathbb{R}}
\newcommand{\p}{\mathcal{P}}
\newcommand{\T}{\mathbb{T}}
\newcommand{\g}{\mathcal{G}}
\DeclareMathOperator{\var}{var}

\def \( {\left( }

\def\){\right)}

\def\com{G\circ\pi}
\def\acom{(\alpha/(\alpha+1))G\circ \pi}

\begin{document}
\title
{Dimensions of compact invariant sets of some expanding maps}
\author{Yuki Yayama}
\address{Centro de Modelamiento Matem\'{a}tico, Universidad de Chile, Av. Blanco Encalada 2120 Piso 7 Santiago de Chile} \email{yyayama@dim.uchile.cl}
\begin{abstract}
We study the Hausdorff dimension and measures of full Hausdorff
dimension for a compact invariant set of an expanding nonconformal
map on the torus given by an integer-valued diagonal matrix.  The
Hausdorff dimension of a ``general Sierpinski carpet" was found by
McMullen and Bedford and the uniqueness of the measure of full
Hausdorff dimension in some cases was proved by Kenyon and Peres. We
extend these results by using compensation functions to study a general Sierpinski carpet
represented by a shift of finite type. We give some conditions under which a general
Sierpinski carpet has a unique measure of full Hausdorff dimension and study the properties of the unique measure.
\end{abstract}
\maketitle

\section{Introduction}
Many natural problems in dynamics are concerned with the following kinds of questions.
Let $T$ be a continuous expanding map on a Riemannian manifold $M,$ and
$K$ be a $T$-invariant compact subset. What is the
Hausdorff dimension of $K$?  Is there an ergodic measure of full
Hausdorff dimension? If so, is it unique? What are the dynamical
properties of any measure(s) of full Hausdorff dimension?
A useful summary of this area may be found in Gatzouras and Peres \cite{GP},
where these questions are explicitly formulated.

The simplest examples are the middle-third Cantor set and the Sierpinski carpet.
These are compact invariant sets for the conformal maps $x\mapsto 3x \textnormal{ mod } 1$ and
$(x,y)\mapsto (3x \textnormal{ mod } 1,$ \\
\noindent $3y \textnormal{ mod } 1),$ respectively (maps that dilate all directions uniformly). The questions above are well understood in the conformal setting (see Gatzouras and Peres \cite{GP} for an overview and history). For the nonconformal setting, only a few cases are understood (see \cite{B, GP, KP, Mc}). In this paper we study the nonconformal case, in particular, consider nonconformal expanding maps of the torus
given by $T(x,y)=(lx \textnormal{ mod }1, my \textnormal{ mod }1 ),
l>m\geq 2, l,m \in \N.$
Bedford \cite{B} and  McMullen \cite{Mc} independently answered the question
on Hausdorff dimension for compact $T$-invariant subsets of the following form.
Consider a subset of the torus constructed in the following way. At
the first step, draw $(l-1)$ vertical lines and $(m-1)$ horizontal
lines in the unit square to get $lm$ congruent rectangles. Shade
some of the rectangles in the square and erase the parts that are
not shaded. At the second step, in each shaded rectangle draw again
$(l-1)$ vertical lines and $(m-1)$ horizontal lines to get  $lm$
congruent rectangles, and shade again the corresponding (smaller)
rectangles as in the first step. Erase parts that are not shaded at
this step. Repeating this process, we get a compact $T$-invariant
set which we call an \em NC carpet \em (for nonconformal). Suppose
$R$ is the set of rectangles chosen to be shaded at the first step.
McMullen \cite{Mc} showed, by finding a Bernoulli measure of full
Hausdorff dimension, that the Hausdorff dimension of the set is
given by
\begin{equation} \label{Hdimension}
\log_m(\sum_{j=0}^{m-1}t^{\log_{l}m}_j),
\end{equation} where $t_j$
is the number of members of $R$ in the $j$' th row in the unit
square. Kenyon and Peres \cite{KP} showed that for an NC carpet
there is a unique measure of full Hausdorff dimension. Using a
coding map constructed by a Markov partition for $T,$ one obtains a
symbolic representation of the carpet which is a full shift on
finitely many symbols.

In this paper, we consider a more general
set whose symbolic representation is a shift of finite type rather than a full shift. We
call such a set an \em SFT-NC carpet. \em  We will give partial
answers to the questions above in this setting. Our results give (under some technical conditions) a formula for the Hausdorff dimension as well as uniqueness and the Bernoulli property for the measure of full Hausdorff dimension.  These positive results narrow the possibilities where an example of nonuniqueness of the measure of full Hausdorff dimension might be found (see Section
\ref{lexamples}).

Using the Ledrappier-Young formula and a Markov partition, Gatzouras
and Peres \cite{GP} translated the problem into one concerning factor maps
$\pi:(X, \sigma_X) \rightarrow (Y,\sigma_Y)$
between symbolic dynamical systems: one seeks measures that
maximize, for a fixed $\alpha >0$ determined by the Lyapunov
exponents, the weighted entropy functional $h_{\mu}(\sigma_{X})
+\alpha h_{\pi \mu}(\sigma_{Y})$.
Shin \cite{S2} showed that if there is a
\em saturated compensation function \em $G\circ\pi$ (see page
\pageref{compensationf} for the definition for the factor map
$\pi$), then for any $\alpha>0$ the set of all shift-invariant measures $\mu$ on $X$ that maximize
$h_{\mu}(\sigma_{X}) +\alpha h_{\pi \mu}(\sigma_{Y})$ is the set of
equilibrium states for the function $\acom.$ A saturated compensation
function helps us make some progress on the problem, giving us a systematic way to approach the problem (see pages \pageref{compensationf}, \pageref{symrep}, and \pageref{HDK}). We will answer
questions about SFT-NC carpets for which a saturated compensation function
exists (when they are represented in symbolic dynamics, see page
\pageref{symrep}). Besides the results of Gatzouras-Peres and Shin,
key ingredients of our approach are the work on grid functions by
Hofbauer \cite{H1} and Markley-Paul \cite{MP}, and on functions in
the Bowen class by Walters \cite{W2005} and on $g$-measures
by Coelho and Quas \cite{CQ}.

Let $X,Y$ be topologically mixing shifts of finite type, and let
$\pi: X\rightarrow Y$ be a factor map, which as usual we may assume
to be a 1-block map (see \cite{LM}).  We begin by considering, in
Section \ref{mainthms}, the case when the alphabet of $Y$ is $\{1,2\}.$
Unless $Y$ is trivial, there is a singleton
clump (some symbol in the alphabet for $Y$ has only one preimage in
the alphabet for $X$), so we assume that $\pi^{-1}\{1\}=\{1\}.$  In Theorem \ref{thm1}, we find a saturated
compensation function which is not necessarily in the Bowen class. In
Theorem \ref{thm2}, we find many situations in which there is a saturated
compensation function in the Bowen class. Examples are provided in
Section \ref{fexamples}. In the first example, a saturated compensation
function is not in the Bowen class but is a \em grid function \em in the sense of Markley and Paul
\cite{MP}. We show that the measure of full dimension is Bernoulli, using a theorem of Coelho and Quas (Theorem \ref{CQg}).

In Section \ref{extensions}, we extend the results of Sections
\ref{mainthms} and \ref{Main2} to the case when the alphabet
of $Y$ is $\{1,\cdots, n\}, n\geq 2$. The approach we use to establish
uniqueness of the equilibrium state is again based on determining
when a saturated compensation function is a grid function or is in
the Bowen class. Let $B$ be the transition matrix of $\pi^{-1}\{2,
3, \cdots, k\}.$  We construct a saturated compensation function and
show the uniqueness of the equilibrium state for it under some
conditions (see Theorem \ref{prop1} for the case $h_\text{top}(X_B)=0$  and Theorem \ref{prop2} for the case $h_\text{top}(X_B)>0$ that
make a compensation function have the form of a grid function). This
approach is generalized in Proposition \ref{hopetobeuseful} to find
more saturated compensation functions that have unique equilibrium
states. Then we study some particular situations that are not
covered by the results mentioned so far, namely the case when the
matrix $B$ has a certain block form (Theorem \ref{prop3}) or a
triangular form (Theorem \ref{arigatou}). Theorem \ref{prop3} treats
the case when the graph of $Y$ has no arrow among symbols other than 1. We see in Theorem \ref{arigatou} that if $B$ is
reducible and all the irreducible components are topologically
mixing (with some additional conditions), then there is a saturated compensation function which is in
the Bowen class.

\section{Background}\label{bg}

We review briefly the setup and previous results that we will need. For notation or terminology not explained here, see \cite{LM, KEP, Wtext}.
Write $\Sigma^{+}_{n}=\{1,2,\cdots, n\}^{\N}$ and $\Sigma_{n}=\{1,2,\cdots, n\}^{\Z}.$ We use standard terminology for symbolic dynamics from \cite {LM}, in particular, $[x_0 \cdots x_{n-1}]$ denotes the cylinder set $\{w:w_0=x_0, \cdots, w_{n-1}=x_{n-1}\}.$
If $X$ is a compact metric space and $T:X \rightarrow X$ is a
homeomorphism, let $M(X,T)$ denote the space of all $T$-invariant
Borel probability measures on $X$. For each $\mu\in M(X,T),
h_{\mu}(T)$ denotes the measure-theoretic \em{entropy }\em of $T$
with respect to $\mu$, and let $h_\text{top}(X)$ be the topological entropy
of $T$. For $f\in C(X)$ and $n\geq 1,$ let
$$(S_{n}f)(x)=\sum_{k=0}^{n-1} f(T^{k}(x)).$$ The \em topological
pressure \em $P$ of $T$ is given by $$P(T,f)=\sup\{h_{\mu}(T)+\int f
d{\mu}\vert \mu\in M(X,T)\},$$ and $\mu \in M(X,T)$ is an \em{equilibrium state} \em for $f$ if $P(T,f)=h_{\mu}(T)+\int f
d{\mu}.$  If $T$ is clear from the context, we write $P_{X}(f)$. Conditions under which a potential function $f$ has a unique
equilibrium state $\mu$ and the properties of the unique equilibrium
states are considered throughout this paper.

Let $(X, \sigma)$ be a one-sided topologically mixing shift of
finite type. A measure $\mu \in M(X,\sigma)$ is a (\em{Bowen}\em) \em {Gibbs measure }\em
corresponding to $f \in C(X)$ if there are constants $C_1, C_2>0$ and $P>0$ such that
\begin{equation*}
C_1\leq \frac{\mu([x_{0}x_{1}\cdots x_{n-1}])}{\exp(-Pn +
(S_{n}f)(x))}\leq C_2
\end{equation*}
for every $ x\in X$ and $n\geq 1$.

For each $f \in C(X)$, the \em{Perron-Frobenius operator, }\em
$L_{f}: C(X) \rightarrow C(X),$ is
defined by $(L_{f}\phi)(x)=\sum_{y\in \sigma^{-1} (x)}e^{f(y)}\phi(y).$

We say that $f \in C(X)$ satisfies the \em{RPF condition} \em if
there are $\lambda>0$, $h\in C(X)$ with $h>0$, and $\nu\in
M(X,\sigma)$ for which $L_{f}h=\lambda h, L_{f}^{*}\nu=\lambda \nu$,
$\int h d \nu =1$, and $$\lim_{m\rightarrow \infty}\vert\vert
\lambda^{-m}L_{f}^{m}\phi-h\int\phi d \nu \vert\vert_{\infty}=0$$ for all
$\phi\in C(X)$. It is well known that if $f$ satisfies the
RPF condition then $\mu= \nu h$ is the unique equilibrium state for
$f$ (see \cite{H1}).

Walters \cite{W2005} introduced a class $Bow (X,\sigma)$ of functions that contains
the functions with summable variation, all of which have unique
equilibrium states. Let $$\var_{n}(f)=\sup\{\vert f(x)-f(y)\vert :x,
y \in X, x_i=y_i \textnormal{ for all } 0\leq i\leq n-1\}.$$ Then
$Bow(X,\sigma)=\{f\in C(X): \sup_{n\geq 1} \var_{n}(S_{n}f)<\infty
\}.$
\begin{theorem} \textnormal{\cite{W2005}}\label{bowen}
If $f\in Bow(X,\sigma)$, then $f$ has a unique equilibrium state
$\mu$, and the natural extension of $(X,\sigma, \mu)$ is
measure-theoretically isomorphic to a Bernoulli system.
\end{theorem}

\begin{theorem}\textnormal{\cite{CQ}}\label{CQg}
Suppose $g$ is a continuous function from $\Sigma_{n}^{+}$
to (0,1) and $\Sigma_{i=1}^{n}g(ix)=1$ for all $x \in \Sigma_{n}^{+}$. For each $i=1,2,\cdots,$ let $a_{i}=(n/2) \var_{i}(g)$.
Suppose that there is $r\geq 1$ such that
$$ \sum_{k=r}^{\infty}\prod_{i=r}^{k}(1-a_i)=\infty.$$
Then there is a unique $g$-measure $\mu_{g}$ corresponding to $g$, and the natural extension of $( \Sigma_{n}^{+},
\sigma, \mu_{g})$ is isomorphic to a Bernoulli system.
\end{theorem}
The definition of $g$-measure can be found in the paper of Coelho and Quas \cite{CQ}. Theorem
\ref{CQg} is valid also for the case of shifts of finite type.\\

Based on the potential functions Hofbauer \cite{H1} considered,
Markley and Paul \cite{MP} extended his ideas, introducing \em grid functions\em. Let $(X,\sigma)$
be a subshift of $(\Sigma_n,\sigma)$. Consider a partition $\p=\{\rho(X), M_{0},M_{1},\cdots
\}$ of $\Sigma_{n}^{+}$ satisfying
\label{grid}
the following conditions: \\

\noindent(a) $\rho(X)=\{x\in \Sigma_{n}^{+}: \textnormal{ there
exists } w\in
X \textnormal{ with } x=w_{0}w_{1}w_{2}\cdots\}.$\\
(b) Each $M_n$ is closed and open for $n=0,1,2 \cdots.$\\
(c) $\lim_{n\rightarrow\infty}M_n=\rho(X)$, where the limit is taken
with respect to the Hausdorff metric on the closed subsets
of $\Sigma_{n}^{+}.$\\
(d) There exists $K_0\in \N$ such that if $B$ is any
cylinder set whose length is larger than
$K_0$ and such that $B\cap\rho(X)=\emptyset$, then there
exists $j(B)$ such that
$B\subset M_{j(B)}.$\\
(e) For each positive integer $s$, there is a length $l_{0}(s)> K_0$
such that $j(B)\geq s$ whenever $B=B'b_{l}=b_{0}b_{1}\cdots
b_{l-1}b_{l}$ is a block with $l\geq l_{0}(s), B\cap
\rho(X)=\emptyset$ and $B'\cap \rho(X) \neq \emptyset$.\\

\noindent Then $$\g(X,\p)=\{g= \sum_{n=0}^{\infty} a_{n} {\chi}_{M_n}
: \lim_{n\rightarrow\infty} a_{n}=0\}\subset C(\Sigma_{n}^{+})$$ is
called the set of  \em{grid functions }\em associated with $(X,
\sigma)$ and the partition $\p.$

\begin{theorem}\textnormal{\cite{MP}} \label{MPkey}
For $g\in \g(X,\p),$ $g$ satisfies the RPF condition if and only if
$P(g)>h_\text{top}(X).$  If $P(g)>h_\text{top}(X),$ then the unique
equilibrium state for $g$ is $\nu h$, where $\nu$ and $h$ are found
by the RPF condition. The function $h$ is constant on any cylinder
set $C$ whose length is larger than $K_{0}$ and such that $C\cap
\rho(X)=\emptyset$.
\end{theorem}
If $(X,\sigma)$ is a one-sided topologically mixing shift of finite
type and $f: X\rightarrow \R$ is a continuous function, then $L^{*}_{f}$ has an eigenmeaure corresponding to a positive eigenvalue. Hence Theorem \ref{MPkey} remains valid in the case of topologically mixing shifts
of finite type.\\

Let $S:X\rightarrow X$, $T: Y\rightarrow Y$ be continuous maps on
compact metric spaces. A map $\pi:(X,S) \rightarrow (Y,T)$ is called
a \em factor map \em if it is a continuous surjection with $\pi\circ
S= T\circ \pi.$ Denote by $P_{X}(S, \cdot)$ and $P_{Y}(T, \cdot)$
the topological pressure functionals of $S$ and $T,$ respectively.
Compensation functions were introduced by Boyle and Tuncel \cite{BT} and
studied by Walters \cite{Wcom} in connection with a relative pressure.
A function \label{compensationf}$F\in C(X)$ is a \em{compensation function }\em
for $(S, T, \pi)$ if $$ P_{X}(S, F+\phi\circ \pi)=P_{Y}(T, \phi)
\textnormal{ for all }\phi\in C(Y).$$ If $F=G\circ\pi \in C(X)$ with
$G\in C(Y)$, then $G\circ \pi$ is a \em{saturated compensation
function}\em.

\begin{theorem}\textnormal{\cite{Wcom}} \label{shincoro}
Let $(X,\sigma_{X})$, $(Y, \sigma_{Y})$ be subshifts, and let $\pi$
be a factor map. For each $n\geq 1$ and $y\in Y$, let $D_n(y)$
consist of one point from each nonempty set $\pi^{-1}(y)\cap
[i_{0}i_{1}\cdots i_{n-1}]$. For $G\in C(Y), G\circ \pi$ is
a compensation function for $\pi$ if and only if
$$
\int_{Y} \limsup_{n\rightarrow\infty}\frac{1}{n}
[\log (e^{(S_nG)(y)}\cdot \vert D_{n}(y)\vert)]d{\nu}=0$$ for all $\nu
\in M(Y,\sigma_{Y}).$
\end{theorem}

The following theorems have key roles in helping to solve the
problems under consideration.
\begin{theorem} \label{LYformula} \textnormal{\textbf{(Special case of the Ledrappier-Young
formula \cite{KP})}} Let $2\leq m<l$, $m,l\in \N$. Let $S$ be defined on the
2-torus $X$ by $S(x,y)=(lx \text{ mod }1, my \text{ mod } 1).$ Let
$\mu$ be an ergodic $S$-invariant measure on $X$. Let $\pi$ be the
projection of $X$ onto the $y$-axis, let $T$ be defined on $Y$ by
$Ty=my \textnormal{ mod }1$, and let $\pi\mu= \mu{\pi}^{-1}.$ Then
\label{LYformula}
$$\dim_{H}\mu= \frac{1}{\log l}h_{\mu}(S)+(\frac{1}{\log m}-\frac{1}{\log l})h_{\pi\mu}(T).$$
\end{theorem}
\begin{theorem}\textnormal{\cite{S2}} \label{shinthm}
Let $(X,\sigma_{X})$, $(Y, \sigma_{Y})$ be subshifts and let $\pi$
be a factor map. Suppose $\pi$ has a saturated compensation function
$G\circ \pi, G \in C(Y).$ For $\alpha \geq 0$, the set of all
shift-invariant measures $\mu$ on $X$ which maximize
$$\phi_{\alpha}(\mu)= h_{\mu}(\sigma_X)+\alpha h_{\pi
\mu}(\sigma_{Y})$$  is the set of equilibrium states of
$(\alpha/(\alpha +1)) G\circ\pi\in C(X)$. In particular, it contains an ergodic measure.
\end{theorem}

We now give the definitions of NC carpets and SFT-NC carpets.  Fix
two positive integers $l$ and $m$, $l>m\geq 2.$ Let $T$ be the
endomorphism of the torus ${\T}^{2}= {\R}^{2}/{\Z}^{2}$ given by
$T(x,y)=(lx \textnormal{ mod }1, my \textnormal{ mod }1).$ Let
$$\p=\{ [\frac{i}{l}, \frac {i+1}{l}]\times[\frac{j}{m},
\frac{j+1}{m}]: 0\leq i \leq l-1, 0\leq j\leq m-1\}$$ be the natural
Markov partition for $T.$ Label $[\frac{i}{l}, \frac
{i+1}{l}]\times[\frac{j}{m}, \frac{j+1}{m}], 0\leq i\leq l-1, 0\leq
j\leq m-1$, by the symbol $(i,j)$. Define $(\Sigma^{+}_{lm},
\sigma)$ to be the full shift on these $lm$ symbols. Consider the coding map $\chi: \Sigma
^{+}_{lm}\rightarrow{T^2},$ defined by \label{codingmap}
\begin{displaymath}
\chi(\{(x_k,
y_k)\}_{k=1}^{\infty})=(\sum_{k=1}^{\infty}\frac{x_k}{l^k},
\sum_{k=1}^ {\infty} \frac{y_k}{m^k}).
\end{displaymath}
\label{symrep} Let $R=\{(a_1,b_1), (a_2,b_2) \cdots, (a_r,b_r)\}$ be a subalphabet of the labels of
$\p.$ The \em NC carpet $K(T,R)$ \em is defined by \label{gscarpet}
\begin{equation*}K(T,R)=\{( \sum_{k=1}^{\infty}\frac{x_k}{l^k},
\sum_{k=1}^{\infty}\frac{y_k}{m^k}): (x_k, y_k)\in R \text{ for all
} k\}.
\end{equation*}
It is a compact $T$-invariant subset of the torus. Denote by $A$ a
transition matrix among the members of $R$,
 so that $A$ is an $r
\times r$ matrix with entries 0 or 1. The \em SFT-NC carpet
$K(T,R,A)$ \em is defined by \label{dsc}
\begin{equation*}
K(T,R,A)=\{(\sum_{k=1}^{\infty}\frac{x_k}{l^{k}},
\sum_{k=1}^{\infty} \frac{y_k}{m^{k}}): (x_k, y_k)\in R, A_{(x_k,y_k)(x_{k+1}, y_{k+1})}=1 \textnormal{ for all } k \}.
\end{equation*}

Now let $(X,\sigma)$ be the subshift on symbols of the members of
$R$ with the transition matrix $A$ as above. Let $\pi$ be the
projection map to the $y$-axis. Let $Y=\pi(X).$ Using Theorem
\ref{LYformula} and Theorem \ref{shinthm}, we get the following.
\begin{corollary}\label{HDK}
Suppose there is a saturated compensation function $G\circ\pi\in
C(X)$ with $G\in C(Y)$ for $\pi:X \rightarrow Y.$ Then the Hausdorff
dimension of the SFT-NC carpet $K(T,R,A)$ is given by
\begin{equation*}
  \dim_{H}K(T, R, A)= \frac{P((\alpha/(\alpha +1)) G\circ\pi)}{\log  m}.
\end{equation*}
Suppose in addition that $\acom$ satisfies the RPF condition. Then
the Hausdorff dimension of the carpet is given by $\log_{m}
\lambda,$ where $\lambda$ is the spectral radius of $L_{\acom}$.
\end{corollary}
\begin{proof}
Let $\mu$ be an ergodic $T$-invariant measure on $K(T,R,A).$ Using
the natural coding map, there is a measure $\bar \mu$ on $X$ mapped
to $\mu.$ Using the Ledrappier-Young formula and the fact that
coding map is bounded to one (see \cite{PT}),
\begin{eqnarray*}
  {\mathrm{dim}}_{H}{ \mu}
   &=& \frac{1}{\log l}h_{\mu}(T) + (\frac{1}{\log m}-\frac{1}{\log l}) h_{\pi \mu} (y\rightarrow my) \\
   &=&\frac{1}{\log l}[h_{\bar \mu}(\sigma _X) + (\log_{m}l-1) h_{\pi \bar \mu}
   (\sigma_Y)].
\end{eqnarray*}
By the proof of Theorem \ref{shinthm} \cite{S2}, for $\alpha>0,$
\begin{equation}\label{key1}
P_{X}(\frac{\alpha}{\alpha +1} G\circ \pi)= \frac{1}{\alpha +1}
\sup _{\mu\in M(X, \sigma_{X})}\{h_{\mu}(\sigma_X)+\alpha h_{\mu \circ \pi
}(\sigma_Y)\}.
\end{equation}
Letting $\alpha=\log_{m}l-1$  and using a theorem of Gatzouras and Peres \cite{GP2}, we get
\begin{equation*}
\begin{split}
\dim_{H}K(T,R,A)&= \sup \{\dim_H \mu : \mu(K(T,R,A))=1, \mu \textnormal{ is $T$-invariant and ergodic.}\}\\
&= \frac{1}{\log l}\sup_{ \mu \in M(X, \sigma_{X})} \{h_{
\mu}(\sigma_{X})+(\log_{m}l-1 )h_{ \pi \mu}(\sigma_{Y})\}\\
&=\frac{1}{\log l}(\alpha +1)P(\frac{\alpha}{\alpha +1}G\circ
\pi)=\frac{P(\acom)}{\log m}.
\end{split}
\end{equation*}
For the second part, we show that $P((\alpha/(\alpha +1))G\circ
\pi)=\log \lambda.$ Let $\varphi= \acom.$  Let $\lambda>0, h\in
C(X)$ be obtained by the RPF condition. Let
$$\bar \varphi=({\alpha}/({\alpha +1})) G\circ \pi +\log h-\log h\circ
\sigma -\log \lambda.$$ Then $\mu$ is the $g$-measure for $g=e^{\bar
\varphi}$ \cite{H1}, and so the same arguments as in the proof of
Corollary 3.3 (i) \cite{Wg} give us the result.
\end{proof}
Specializing to NC carpets, we have an alternative, dynamical proof
of the formula given by McMullen \cite{Mc}: The Hausdorff dimension for each NC carpet $K(T,R)$ is given by
Formula (\ref{Hdimension}).

\section{Main Result-Part 1}\label{mainthms}
To formulate and prove our main results, Theorem \ref{thm1} and
Theorem \ref{thm2}, we will employ the following setting (Setting
(A)) and Convention, and we will need the following additional hypothesis [C].\\

\noindent \textbf{\underline{Setting (A)}} Fix $r=3,4,\cdots.$ Let
$X \subset \{1,2,\cdots, r\}^{\N}$ be a topologically mixing shift
of finite type with positive entropy, $Y= \{1,2\}^{\N}$ or $Y\subset
\{1,2\}^{{\N}}$ a shift of finite type with positive entropy, and
$\pi : X \rightarrow Y$ a one-block factor map such that $\pi^{-1}
\{1\}=\{1\}.$ Let $A$ be the transition matrix of $X$ and let $B$ be
the $(r-1)\times (r-1)$ submatrix of $A$ corresponding to the
indices $2,3,\cdots, r$ (giving the transitions among the symbols in
$\pi^{-1} \{2\}).$ Denote by $X_B$ the shift of
finite type determined by $B$. Let $0\leq \tau<1.$\\

\noindent \textbf{Convention:} \label{convch3}For a block $y_{0}\cdots y_{k-1}$ of
length $k$ on $\{1,2\}^{\N},$ we define $\vert \pi^{-1}[y_{0}\cdots y_{k-1}]\vert$ to be
the number of blocks of length $k$ in $X$ mapped to $y_{0}\cdots y_{k-1}$ under $\pi.$
We define $\vert \pi^{-1}[12^{n}1]\vert =1$ if $\pi^{-1}[12^{n}1]
={\emptyset}.$\\

\noindent \textbf{Condition [C]:  } \label{conditionC}Let $\overline
M_n= \pi^{-1}[2^{n}1](=\{x\in X:\pi(x_0 \cdots x_{n})=2^{n}1\}) \text { for } n\geq 1$ and $\overline
Z=\pi^{-1}\{2^{\infty}\}.$
We assume the following:\\
(1) $\overline Z$ is a one-sided shift of finite type such that if
$T$ is the transition matrix of $\overline Z$, $\Sigma_{T}$ is the
two-sided shift of finite type on $\{2, \cdots, r\}^{\Z}$ determined
by $T$, and $\Sigma_{T}\vert _{+}$ is the projection of $\Sigma_{T}$
to a one-sided subshift, then $\Sigma_{T}\vert _{+}= \overline
Z$.\\
(2) $\lim_{n \rightarrow \infty}\overline M_n=\overline Z,$ where
the limit is
taken in the Hausdorff metric.\\

\begin{theorem}\label{thm1}

\renewcommand{\theenumii}{\alph{enumii}}
\renewcommand{\theenumi}{\arabic{enumi}}
Let Setting (A) hold, and suppose that $h_\text{top}(X_{B})=0.$
     \begin{enumerate}
     \item Suppose that there is $n\geq 1$ such that $B^n=0$. Then there is a compensation function which is
     locally constant.
     \item Suppose for every $n$ that $B^n\neq 0$ and that the following
     two conditions hold:\\
        (i)
       \begin{equation*}\label{limquotient1}
       \lim_{n \to\infty} \frac{\vert \pi^{-1}[12^{n-1}1]\vert}{\vert
       \pi^{-1}[12^{n}1]\vert}=1.
       \end{equation*}
       (ii)
       \begin{equation*}\label{-topofXB1}
        h_\text{top}(X_B)=\lim_{n \to\infty}\frac{\log \vert \pi^{-1}[2^{n}]\vert}{n}=0.
       \end{equation*}

\noindent Then
            \begin{enumerate}
            \item  There exists a compensation function $G\circ \pi \in C(X)$
            such that $G \in C(Y)$.
            \item  $\tau G\circ \pi$ has a unique equilibrium state $\mu$.
            \item Under Condition [C] above, $(\sigma, \mu)$ is exact, hence strongly mixing.
            \item The unique equilibrium state $\mu$ is Gibbs if and only if
            $\sup_{n} \vert \pi^{-1}[12^{n}1]\vert < \infty$.
            \item If the
            unique equilibrium state $\mu$ is Gibbs, then the natural extension of
            $(\sigma,\mu)$ is isomorphic to a Bernoulli system.
            \end{enumerate}
\end{enumerate}
\end{theorem}

\begin{remark}
\textnormal{The hypothesis that $X$ be topologically mixing is not necessary
for (2)(a).}
\end{remark}

\subsection {Proof of (2)(a).} \label{pfmain12a} \noindent Define $G:
Y\rightarrow \R$ by
\begin{equation}\label{comfun} G(y)= \begin{cases}
 \log (\vert \pi^{-1}[12^{k-1}1]\vert/ \vert
\pi^{-1}[12^{k}1]\vert) & \textrm{
 if }y\in[2^{k}1], k\geq 2 \\
\log ({1}/{\vert \pi^{-1}[121]\vert}) & \textrm{
 if }y\in[21] \\
0 & \textrm{
 if }y\in [1] \textnormal{ or } y=2^{\infty}.
\end{cases}
\end{equation}
We show that $G\circ \pi$ is a compensation function for $\pi$ by
showing that
\begin{equation}\label{cosc}
\int_{Y} \limsup_{n\rightarrow \infty} \frac{1}{n}[\log
(e^{(SnG)(y)}\vert D_n(y)\vert )] d\nu= 0 \textnormal { for all }
\nu \in M(Y,\sigma_{Y})
\end{equation}
(see Therorem \ref{shincoro} for the definition of $D_n$). For
$i=1,2,$ let
$$E_i=\{y \in Y: y=y_0 \cdots y_{p-1} i^{\infty} \textnormal{ for
some }p\geq 1, y\neq i^{\infty}\}.$$ Then for any $\nu\in M(Y,\sigma_{Y})$,
$\nu(E_1)=\nu(E_2)=0$ because $\nu$ is a $\sigma_{Y}$-invariant measure.
It is enough to show that
\begin{equation}\label{Nice eq}
\limsup_{n\rightarrow \infty} \frac{1}{n} [\log (e^{(SnG)(y)}\vert
D_n(y)\vert )]= 0 \textnormal { for all } y \in Y \setminus (E_1
\cup E_2).
\end{equation}

\subsubsection{Case 1-1.} \label{maincase1-0} We consider the case when
$y\notin E_1 \cup E_2 \cup \{1^{\infty},2^{\infty}\}.$ Let $n>2$ be
fixed and consider the first $n$ states $y_0y_1\cdots y_{n-1}$ of
$y\in Y.$ Let $y_0=2$ and $y_{n-1}=1$. Fix such a $y\in Y.$ Then
there exist $k_1 ,k'_1, \cdots ,k_l ,k'_l\geq 1 ,$ where $k_i$ for
$2\leq i \leq l$, $k'_i$ for $1\leq i \leq l,$ and $l$ depend on $y$
and $n,$ and $k_1$ depends on $y$, such that $k_1 +k'_1+ \cdots +k_l
+k'_l =n,$ and $y= 2^{k_1}1^{k'_1}\cdots 2^{k_l}1^{k'_l}\cdots.$
Then
$$\vert
D_n(y)\vert = \vert \pi^{-1}[2^{k_1}1]\vert \vert
\pi^{-1}[12^{k_2}1]\vert\cdots \vert \pi^{-1}[12^{k_l}1]\vert$$ and
\begin{equation*}
(S_{n}G)(y)=\log \frac{1}{\vert \pi^{-1}[12^{k_1}1]\vert \vert
\pi^{-1}[12^{k_2}1]\vert \cdots \vert \pi^{-1}[12^{k_l}1]\vert}.
\end{equation*}
\label{sng} Therefore,
$$
\vert D_n(y)\vert e^{(S_n G)(y)}= \frac{\vert
\pi^{-1}[2^{k_1}1]\vert}{\vert \pi^{-1}[12^{k_1}1]\vert}, $$ which
is a constant depending on $y$. Hence (\ref{Nice eq}) holds for such
$y$.

\subsubsection{Case 1-2.}\label{maincase1-1} Let $n$ be fixed as in \ref{maincase1-0}.
Consider $y\notin E_1 \cup E_2 \cup \{1^{\infty},2^{\infty}\}$ and
$y_0=y_{n-1}=2.$
 Fix $y\in Y.$ Then there exist $k_1
,k'_1, \cdots ,k_{l-1},k'_{l-1},k_l \geq 1$ (where $k_i$ for $1\leq
i \leq l$, $k'_i$ for $1\leq i \leq l-1$, and $l$ depend on $y$ and
$n$), and there exists $t\geq 0$ (where $t$ depends on $y$ and $n$),
such that $k_1 +k'_1+ \cdots +k_{l-1} +k'_{l-1}+k_l =n$ and $y=
2^{k_1}1^{k'_1}\cdots 2^{k_{l-1}}1^{k'_{l-1}}2^{k_l}2^{t}1\cdots.$
Depending on the values of $k_l$ and $t$ associated with $y$ and
$n$, $n$ falls into one of four subsequences of $\N$ (see below). We
will show that the estimate for the limsup holds along each of these
four subsequences of $n$'s, and hence holds along the full sequence.

We have
$$\vert D_n(y)\vert \leq  \vert \pi^{-1}[2^{k_1}1]\vert \vert \pi^{-1}[12^{k_2}1]\vert
\cdots \vert \pi^{-1}[12^{k_{l-1}}1]\vert \vert
\pi^{-1}[12^{k_{l}}]\vert $$ and

\begin{equation}\label{Sn22}
\begin{split}
(S_{n}G)(y) = \log \frac{\vert \pi^{-1}[12^{t}1]\vert }{\vert
\pi^{-1}[12^{k_1}1]\vert
 \vert \pi^{-1}[12^{k_2}1]\vert\cdots
\vert \pi^{-1}[12^{k_{l-1}}1]\vert  \vert
\pi^{-1}[12^{k_{l}+t}1]\vert}.
\end{split}
\end{equation}
\label{sng12}
 If we let $C(y)= {\vert
\pi^{-1}[2^{k_1}1]\vert}/{\vert \pi^{-1}[12^{k_1}1]\vert},$ then

\begin{equation}\label{ineq1}
\vert D_n(y)\vert e^{(SnG)(y)} \leq C(y)\frac{\vert
\pi^{-1}[12^{t}1]\vert \vert \pi^{-1}[12^{k_l}]\vert}{\vert
\pi^{-1}[12^{k_{l} +t}1]\vert}.
\end{equation}

Now fix $\epsilon>0$.  Then by (2)(i), there exists $N \in \N$ such
that
$$ (1-\epsilon)\vert\pi^{-1}[12^{t+1}1]\vert \leq
\vert\pi^{-1}[12^{t}1]\vert\leq
(1+\epsilon)\vert\pi^{-1}[12^{t+1}1]\vert $$ for all $t\geq N.$
Therefore, in the case when $t=t(y,n)\geq N,$ $k_{l}=k_{l}(y,n)\geq
1$,
$$\vert\pi^{-1}[12^{t}1]\vert\leq
(1+\epsilon)^{k_{l}}\vert\pi^{-1}[12^{t+k_l}1]\vert,$$ and so
\begin{equation} \label{approx1}
\frac
{\vert\pi^{-1}[12^{t}1]\vert}{\vert\pi^{-1}[12^{t+k_l}1]\vert}\leq
(1+\epsilon)^{k_{l}} \leq (1+\epsilon)^{n}.
\end{equation} For the
same fixed $\epsilon>0$, by (2)(ii), there exists $\overline N \in
\N$ such that
\begin{equation}\label{klarge}
\vert \pi^{-1}[2^{k}]\vert\leq e^{k \epsilon} \textnormal{ for all }
k \geq \overline N .
\end{equation}
For $1\leq k \leq \overline N -1$, if we let $\beta= \max _{1 \leq k
\leq {\overline N -1}}  (\log \vert \pi^{-1}[2^{k}]\vert)/{k},$ then
\begin{equation}\label{ksmall}
\vert \pi^{-1}[2^{k}]\vert \leq e^ {\beta \overline N}\textnormal{
for all }1\leq k \leq {\overline N -1}.
\end{equation}
Now consider $n$ such that $t=t(y,n)\geq N$ and $k_l=k_{l}(y,n) \geq
\overline N.$  By (\ref{approx1}), (\ref{klarge}), $\vert
\pi^{-1}[12^{k_l}]\vert \leq \vert \pi^{-1}[2^{k_l}]\vert$, and
$k_l\leq n,$ we have for such $n$ that

\begin{equation*}
\vert D_n(y)\vert  e^{(SnG)(y)} \leq  C(y) (1+\epsilon)^n \vert
\pi^{-1}[12^{k_{l}}]\vert \leq C(y)(1+\epsilon)^n {e^{n\epsilon}}.
\end{equation*}
\label{sequence1} Taking the increasing subsequence
$\{n_i\}^{\infty}_{i=1}$ of such $n$, we have
\begin{equation}\label{tk1}
\begin{split}
\limsup_{i\rightarrow \infty } \frac {1}{n_i} \log
(e^{(Sn_{i}G)(y)}\vert D_{n_{i}}(y)\vert ) &\leq
\limsup_{i\rightarrow \infty }\frac {
n_{i} \log(1+\epsilon)+n_{i} \epsilon}{n_{i}}\\
& = \log(1+\epsilon)+\epsilon < 2\epsilon.
\end{split}
\end{equation}
\label{sequence2} Now consider $n$ such that $t=t(y,n)\geq N$ and
$1\leq k_l=k_l(y,n)<\overline N.$ By (\ref{approx1}),
(\ref{ksmall}),$\vert \pi^{-1}[12^{k_l}]\vert \leq \vert
\pi^{-1}[2^{k_l}]\vert $, and $k_l< \overline N,$ we have for
such $n$ that
\begin{equation*}
 \vert D_n(y)\vert  e^{(SnG)(y)} \leq  C(y)
(1+\epsilon)^n \vert \pi^{-1}[12^{k_{l}}]\vert \leq
C(y)(1+\epsilon)^{n} e^{\beta \overline N}.
\end{equation*}
Taking the increasing sequence $\{n_i\}^{\infty}_{i=1}$ of such $n$,
we have that
\begin{equation}\label{tk2}
\begin{split}
\limsup_{i\rightarrow \infty } \frac {1}{n_i} \log
(e^{(Sn_{i}G)(y)}\vert D_{n_{i}}(y)\vert ) &\leq
\limsup_{i\rightarrow \infty }\frac {
n_i \log(1+\epsilon)}{n_i}\\
& = \log(1+\epsilon)<\epsilon.
\end{split}
\end{equation}
Now consider $n$ such that $ 0\leq t=t(y,n) < N$ and $
k_l=k_l(y,n)\geq \overline N.$ Let $M=\max_{0\leq t\leq N-1}\vert
\pi^{-1}[12^{t}1]\vert.$ Then by using (\ref{klarge}) and $1\leq
\vert\pi^{-1}[12^{t+k_l}1]\vert,$ referring to (\ref{ineq1}), we get
for such $n$ that
\begin{equation*}
\vert D_n(y)\vert  e^{(SnG)(y)} \leq  C(y)M\vert
\pi^{-1}[2^{k_l}]\vert \leq  C(y)M e^{n \epsilon}.
\end{equation*}
Thus taking the increasing subsequence $\{n_i\}^{\infty}_{i=1}$ of
such $n$,
\begin{equation}\label{tk3}
\limsup_{i\rightarrow \infty } \frac {1}{n_i} \log
(e^{(Sn_{i}G)(y)}\vert D_{n_{i}}(y)\vert ) \leq \epsilon.
\end{equation}
Consider $n$ such that $0\leq t=t(y,n)< N$ and $1\leq
k_l=k_{l}(y,n)<\overline N.$  By (\ref{ksmall}), for such $n$, we
have
\begin{equation*}
\vert D_n(y)\vert  e^{(SnG)(y)} \leq  C(y)M\vert
\pi^{-1}[2^{k_l}]\vert \leq  C(y)M e^{\beta \overline N}.
\end{equation*}
Thus for the increasing subsequence $\{n_i\}^{\infty}_{i=1}$ of such
$n$,
\begin{equation}\label{tk4}
\limsup_{i\rightarrow \infty } \frac {1}{n_i} \log
(e^{(Sn_{i}G)(y)}\vert D_{n_{i}}(y)\vert ) \leq 0.
\end{equation}
By (\ref {tk1}) through (\ref {tk4}), for $y\notin E_1 \cup E_2 \cup
\{1^{\infty},2^{\infty}\}$ and $y_0=y_{n-1}=2,$
$$\limsup_{n\rightarrow \infty } (1/n)[\log (e^{(SnG)(y)}\vert D_n(y)\vert )]<
2\epsilon. $$ By \em Case 1-1 \em  and \em Case 1-2\em, for  $y_0=2$
and $y\notin E_1 \cup E_2 \cup \{1^{\infty},2^{\infty}\},$   we get
$$\limsup_{n\rightarrow \infty} (1/n)[\log(e^{(SnG)(y)} \vert
D_n(y)\vert )]\leq 0.$$
Note that from above in \em Case 1-1 \em we
have that $\vert D_n(y)\vert e^{(SnG)(y)}\geq 1 \textnormal { for }
y \textnormal{ with } y_{n-1}=1.$ \label{case120}By the definition
of limsup, for $y$ such that $y_0=2$ and $y\notin E_1 \cup E_2 \cup
\{1^{\infty},2^{\infty}\},$ we also have
$$\limsup_{n\rightarrow \infty}
(1/n)[\log(e^{(SnG)(y)} \vert D_n(y)\vert )] \geq 0.$$ Therefore, we
get (\ref{Nice eq}) for all $y \notin E_1 \cup E_2
\cup\{1^{\infty}, 2^{\infty}\}$ and $y_0=2.$\\

\subsubsection{ Case 1-3.} Let $n>2$ be fixed, as in \ref{maincase1-0}. Consider $y\notin
E_1 \cup E_2 \cup \{1^{\infty},2^{\infty}\}$ and $y_0=y_{n-1}=1$.
Let $y= 1^{k}2^{k_1}1^{k'_1}\cdots 2^{k_l}1^{k'_l}\cdots$ as in
\ref{maincase1-0}, where $k+k_1 +k'_1+\cdots +k_l +k'_l =n.$ Then we
get $\vert D_n(y)\vert e^{(SnG)(y)}=1.$

\subsubsection{ Case 1-4. } \label{maincase1-3}Let $n>2$ be fixed, as in \ref{maincase1-0}. Consider
$y\notin E_1 \cup E_2 \cup \{1^{\infty},2^{\infty}\}$ and $y_0=1$
and $y_{n-1}=2.$  Let $y= 1^{k}2^{k_1}1^{k'_1}\cdots
1^{k'_{l-1}}2^{k_l}2^{t}1\cdots,$ as in \ref{maincase1-1}, where
$k+k_1 +k'_1+\cdots +k'_{l-1}+k_l =n.$ Then we get
$$\vert D_n(y)\vert e^{(SnG)(y)}\leq  \frac {\vert
\pi^{-1}[12^{t}1]\vert \vert \pi^{-1}[12^{k_l}]\vert}{\vert
\pi^{-1}[12^{k_{l} +t}1]\vert} \textnormal{ ($C(y)$ is replaced by
1)}.$$

By \em Case 1-3 \em and \em Case 1-4\em, we have (\ref{Nice eq}) for
for all $y\notin E_1\cup E_2 \cup \{1^{\infty},2^{\infty}\},$ and
$y_0=1.$

\subsubsection{ Case 2.}\label{nicecase} We consider the case when $y\in \{1^{\infty},
2^{\infty}\}.$ If $y=2^{\infty},$ we use $1\leq \vert D_n(y)\vert \leq \vert
\pi^{-1}[2^{n}]\vert $ and $(S_{n} G)(y)=0$ for all
$n.$ For $y=1^{\infty},$ we use $\vert D_n(y)\vert$=1 and
$(S_{n}G)(y)=0$.

Therefore, by \ref{maincase1-0} through \ref{nicecase} , (\ref{Nice
eq}) is satisfied for all $y\in Y \setminus (E_1 \cup E_2).$

\subsection {Proof of (1).} Define $G$ as in
(\ref{comfun}). Since $B^n=0$ for some $n\geq 1,$ there exists
$n_1>0$ such that, using our convention, $ \vert
\pi^{-1}[12^{n}1]\vert=1$ for all $n\geq n_{1}.$ Thus $G$ is locally
constant and  (\ref{Nice eq}) holds for $y\in Y \setminus (E_1 \cup
E_2 \cup{\{2^{\infty}\}}).$

\begin{remark}
\textnormal{Shin \cite{Shin1} gave a sufficient condition for the existence of a
saturated compensation function for the case when
$X\subset\{1,2,\cdots, r\}^{\Z}$ is a topologically mixing shift of
finite type, $Y=\{1,2\}^{\Z}$ or $Y$ is the golden mean subshift,
and $\pi^{-1}\{1\}=\{1\}.$  The compensation function is defined on
the \em two-sided \em shift space $Y.$ The construction of a
compensation function on \em one-sided \em subshifts in this paper
is slightly different, and our sufficient conditions for the
existence of a saturated compensation functions are also slightly
different. But the approach here to show the existence of a
saturated compensation function by using Theorem \ref{shincoro}
is the same as the approach in \cite{Shin1}, and there may be
relations between these two compensation functions.}
\end{remark}
\subsection{Proof of (2)(b).} \label{pfu0} Let
$Z$=$\{2^{\infty}\}$. Define $M_0= [1]$ and $M_n=[2^n 1]$ for all
$n\geq 1.$ Then $\rho (Z)=\{2^{\infty}\}.$ Consider the partition
$\p=\{\rho (Z), M_0, M_1\cdots\}$ of $Y$. Then, using hypothesis
(2)(i), $G\in \g(Z,\p)$ (see page \pageref{grid}). We will use three
Theorems from \cite{PQS}, \cite{Shin1}, and \cite{Wcom}, and Formula
(\ref{key1}) (in the proof of Theorem \ref{shinthm} \cite{S2}) which
yield a corollary that helps to prove the uniqueness of the
equilibrium state for $\acom$.

\begin{theorem}\textnormal{\cite{Wcom}} \label{CC}
Suppose $(X,\sigma_{X})$, $(Y,\sigma_{Y})$ are subshifts and $\pi:
(X,\sigma_{X}) \rightarrow (Y,\sigma_Y)$ is a factor map. Suppose
$F\in C(X)$ is a compensation function. Let $\mu\in M(X, \sigma_{X})$ and
$\phi\in C(Y).$ Then $\mu$ is an equilibrium state of $F+\phi\circ
\pi$ if and only if $\pi\mu$ is an equilibrium state of
$\phi$ and $\mu$ is a relative equilibrium state of $F$ over
$\pi\mu.$
\end{theorem}

\begin{theorem} \textnormal{\cite{PQS}} \label{PQS1}
Let $(X,\sigma_{X})$ be an irreducible shift of finite type,
$(Y,\sigma_{Y})$ a subshift, and $\pi: (X,\sigma_{X})\rightarrow
(Y,\sigma_{Y})$ a factor map. Suppose that there is a symbol $a$ of
$Y$ whose inverse image is a singleton, which is also denoted by
$a$. Then every ergodic measure on $Y$ which assigns positive
measure to [$a$] has a unique preimage of maximal relative entropy.
\end{theorem}

\begin{theorem}\textnormal{\cite{Shin1}}\label{G(2)}
Let $X \subset \{1,2,\cdots, r\}^{\Z}$ be a shift of finite type,
$Y= \{1,2\}^{\Z}$ a full shift or $Y\subset \{1,2\}^{{\Z}}$ the golden mean
shift, and $\pi : X \rightarrow Y$ a factor map such that
$\pi^{-1} \{1\}=\{1\}.$ Suppose $G\in C(Y)$ and $G\circ\pi$ is a
compensation function. Then $G(2^{\infty})=-\limsup_{n} ((\log \vert
\pi^{-1}[2^n]\vert)/n).$
\end{theorem}
Note that Theorem \ref{G(2)} is valid when $X$ is a one-sided
shift of finite type.
\begin{corollary}\label{unique}
Fix $r=3,4,\cdots.$ Let $X \subset \{1,2,\cdots, r\}^{\N}$ be a
topologically mixing shift of finite type, $Y= \{1,2\}^{\N}$ or
$Y\subset \{1,2\}^{{\N}}$ a shift of finite type, and $\pi : X
\rightarrow Y$ a one-block factor map such that $\pi^{-1}
\{1\}=\{1\}.$ Suppose there is a saturated compensation function
$G\circ \pi$ for $\pi.$  Then, for any  $\alpha>0,$ if
$-({1}/({\alpha +1}))G$ has a unique equilibrium state, then
$(\alpha/(\alpha +1))G\circ \pi$ has a unique equilibrium state.
\end{corollary}
\begin{proof}
Let $\nu$ be an ergodic equilibrium state for $-({1}/(\alpha+1))G.$
Let $\mu$ be a preimage of $\nu$ with maximal relative entropy. We
first show that $\mu$ is a relative equilibrium state of $G\circ
\pi$ over $\nu.$ Using the relative variational principle
\cite{Wcom}, we have
\begin{equation*}
\begin{split}
\int P(\sigma_{X}, \pi, G\circ \pi) d\nu &= \sup
\{h_{\mu}(\sigma_{X} \vert _{\sigma_{Y}})+ \int G\circ \pi d\mu:
\mu\in M(X,\sigma_{X}), \pi\mu=\nu\}\\
&=\sup \{h_{\mu}(\sigma_{X})-h_{\pi\mu}(\sigma_{Y})+
\int G\circ \pi d\pi\mu: \mu \in M(X,\sigma_{X}),\\
& \quad \pi\mu=\nu\}\\
&=\sup \{h_{\mu}(\sigma_{X}): \mu \in M(X,\sigma_{X}), \pi\mu =\nu\}-h_{\nu}(\sigma_{Y}) +\int G d\nu.
\end{split}
\end{equation*}
Therefore, $\mu$ is a preimage of maximal entropy if and only if it
is a relative equilibrium state of $G\circ \pi$ over $\nu.$  By
Theorem \ref{CC}, using $\phi= -({1}/({\alpha+1}))G$ and $F=G\circ
\pi$, $\mu$ is an equilibrium state of $(\alpha/({\alpha +1}))
G\circ \pi.$\label{preimagekey}

Next we show that $\nu([1])>0. $ Assume $\nu([1])=0.$ Since
$\nu(2^{\infty})=1,$
\begin{equation*}
\begin{split}
 P_{Y}(-\frac{1}{\alpha +1}G)&= h_{\nu}(\sigma_{Y})-\int
\frac{1}{\alpha +1}G d\nu = -\frac{1}{\alpha +1} G(2^{\infty})\\
&= \frac{1}{\alpha +1} h_\text{top}(X_B) \text{  (by Theorem
\ref{G(2)})}.
\end{split}
\end{equation*}
Take the Shannon-Parry measure ${\mu}_{\max}.$ Since $X_B\subsetneq
X$ and $X$ is topologically mixing, we have
$h_\text{top}(X)>h_\text{top}(X_B).$ Using Formula (\ref{key1}), page
\pageref{key1}, and the definition of compensation function,
\begin{equation}\label{key2}
\begin{split}
P_{Y}(-\frac {1}{\alpha +1}G)&=P_{X}(\frac{\alpha}{\alpha+1}G\circ
\pi) \\
& \geq  \frac {1}{\alpha+1} h_{{\mu}_{\max}}(\sigma_{X}) =\frac
{1}{\alpha+1} h_\text{top}(X) >\frac {1}{\alpha+1} h_\text{top}(X_B).
\end{split}
\end{equation}
This is a contradiction. Now $\nu$ has a unique preimage of maximal
relative entropy by Theorem \ref{PQS1}.

By hypothesis $\nu$ is the unique equilibrium state for
$-({1}/(\alpha +1))G$. Suppose now that there are two equilibrium
states $\mu_1, \mu_2$ for $(\alpha/(\alpha +1)) G\circ \pi.$ Then,
by the above, they are distinct preimages of maximal relative
entropy of $\nu$, and this is a contradiction. \end{proof} Now we
finish the proof of (2)(b). By Corollary \ref{unique}, it is enough
to show that $-({1}/({\alpha +1}))G$ has a unique equilibrium state.
Since $-({1}/({\alpha +1}))G \in \g(Z, \p),$ using (\ref{key2}) and
$h_\text{top}(Z)=0,$
$$P_{Y}(-\frac {1}{\alpha +1}G)>h_\text{top}(Z).$$
By Theorem \ref{MPkey},  $-({1}/({\alpha +1}))G$ has a unique
equilibrium state.
\begin{remark} \textnormal{If $Y$=$\Sigma^{+}_{2}$ and $G$ is as above, then
$tG$ has a unique equilibrium state for any $t \in \R$ by \cite{H1}.}
\end{remark}
\subsection{Proof of (2)(d) and (e).} \label{WH}
\begin{lemma}\label{BandH} Suppose $(X,\sigma)$ is a one-sided topologically mixing shift of finite type. Let
$\varphi \in C(X).$ Then an equilibrium state for $\varphi$ is Gibbs
if and only if $\varphi \in Bow(X, \sigma)$. In particular, if an
equilibrium state for $\varphi$ is Gibbs, then it is the unique
equilibrium state.
\end{lemma}

\begin{proof}
\noindent We follow the arguments in \cite{H1}.  Suppose there is an
equilibrium state $\mu$ for $\varphi$ which is Gibbs. Then there
exist $\lambda, C_1, C_2>0$ such that
$$ C_1\leq \frac{\mu([x_0x_1\cdots x_{n-1}])}{\lambda^{-n} e^{(S_{n}
\varphi)(x)}}\leq C_2 $$ for every $x \in X, n \in \N$.  For $x,x' \in
 [x_0x_1\cdots x_{n-1}]$, we have
$$C_1 e^{(S_{n} \varphi)(x)}\leq \lambda^{n} \mu([x_0x_1\cdots
x_{n-1}])\leq C_2 e^{(S_{n} \varphi)(x')}.$$ Taking logarithms,  we
get
$$ \vert (S_{n} \varphi)(x)-(S_{n} \varphi)(x') \vert \leq \log
\frac{C_2}{C_1}.$$ Therefore, $ \sup_{n\geq 1} \var_{n}(S_{n}\varphi) \leq
\log (C_2/C_1)$ and so $\varphi \in Bow(X, \sigma)$.

Conversely, suppose that $\varphi \in Bow(X, \sigma)$.  Let
$\lambda= e^{P_{X}(\varphi)}.$ Then there exists a unique measure $\nu
\in M(X,\sigma)$ with $L^{*}_{\varphi} \nu= \lambda \nu$, and the unique
equilibrium state for $\varphi$ is $\mu= \nu h$, where
$h:X\rightarrow [a, b]$ is a measurable function, $0<a<b$,
$L_{\varphi} h= \lambda h$ and $\int h d \nu =1$ \cite{W2005}. Since
the shift of finite type $X$ with transition matrix $A$ is
topologically mixing and of positive entropy, there is a number $k$
such that $ A^{k}>0$. Let $r$ be the number of symbols in the
alphabet for $X$. Then, by \cite{W2005}, $\nu$ satisfies
\begin{equation*}\label{gibbseq}
\frac {1}{e^{\var_n( S_{n} \varphi)} \lambda^{k} e^{k \| \varphi \|
}}\leq \frac{\nu([x_0x_1\cdots x_{n-1}])}{\lambda^{-n}
e^{(S_{n} \varphi)(x)}}\leq r^{k}e^{\var_n( S_{n}
\varphi)}\frac{e^{k \|\varphi \|}}{\lambda ^{k}}.
\end{equation*}  Since $\varphi \in
Bow(X, \sigma)$, we have $\var_n (S_{n}\varphi)\leq L$ for all $n$
for some $L$. Therefore $\nu$ is Gibbs and so $\mu=\nu h$ is Gibbs.
\end{proof}
\label{notgibbs}
\begin{lemma}\label{claim1}
Under the assumptions of (2)(i) and (ii), $\vert
\pi^{-1}[12^{n}1]\vert$ is bounded for all $n$ if and only if the
unique equilibrium state of $\varphi= \tau G\circ \pi$ is Gibbs.
\end{lemma}
\begin{proof}
\noindent Suppose first that $\vert \pi^{-1}[12^{n}1]\vert$ is not
bounded. Assume that the unique measure $\mu$ is Gibbs. Consider
$y\in X$ such that $\pi(y)\in[2^{n}1]$ and $z\in X$ such that $\pi(z)=2^{\infty}$. Then $\pi(y),\pi(z)\in[2^{n}]$
and $\vert (S_{n} \varphi)(y)-(S_{n} \varphi)(z) \vert = \tau \log {\vert
\pi^{-1}[12^{n}1]\vert}$, referring to (\ref{Sn22}) and the fact
that $(S_n \varphi)(z)=0$ for all $n$, and $\var_{n}(S_n
\varphi)$ is not bounded. This is a contradiction.

Conversely, suppose that $\vert \pi^{-1}[12^{n}1]\vert \leq K$ for
all $n$. We show  $\varphi=\tau G\circ \pi \in Bow(X, \sigma)$, that
is, $\sup_{n\geq 1}\var_{n}(S_{n} \varphi)< \infty.$

\em {Case 1.} \em Consider $x$ and $x'$ such that
$\pi(x_i)=\pi(x'_i)$ for all $0\leq i\leq n-1$ and
$\pi(x_{n-1})=\pi(x'_{n-1})=2.$

(i) Suppose $\pi(x)=1^{k_1}2^{l_1}1^{k_2}2^{l_2}\cdots
1^{k_t}2^{\beta}1\cdots$ , where $k_1 +l_1+\cdots +l_{t-1}+k_{t}=k<n,
k+\beta \geq n, \pi(x_{n-1})=2$, and
$\pi(x')=1^{k_1}2^{l_1}1^{k_2}2^{l_2}\cdots
1^{k_t}2^{\gamma}1\cdots$ , where $k_1 +l_1+\cdots +l_{t-1}+k_{t}=k<n,
k+\gamma \geq n, \pi(x'_{n-1})=2.$ Direct computation shows that
\begin{equation}\label{eq1}
\vert (S_{n} \varphi)(x)-(S_{n} \varphi)(x') \vert = \tau\vert \log
\frac{ \vert\pi^{-1}[12^{\beta-(n-k)}1]\vert
\vert\pi^{-1}[12^{\gamma}1]\vert} {\vert\pi^{-1}[12^{\beta}1]\vert
\vert\pi^{-1}[12^{\gamma-(n-k)}1]\vert}\vert.
\end{equation}
\label{vratio}
Since $1\leq \vert \pi^{-1}[12^{n}1]\vert \leq K$ for
all $n$, clearly, for any $\beta, \gamma, n,k$ above
\begin{equation}\label{eq1.1}
0\leq \vert (S_{n} \varphi)(x)-(S_{n} \varphi)(x') \vert\leq 2\tau
\log K.
\end{equation}

(ii) Suppose $\pi(x)=1^{k_1}2^{l_1}1^{k_2}2^{l_2}\cdots
1^{k_t}2^{\beta}1\cdots$ , where $k_1 +l_1+\cdots +l_{t-1}+k_{t}=k<n,
k+\beta \geq n, \pi(x_{n-1})=2$, and
$\pi(x')=1^{k_1}2^{l_1}1^{k_2}2^{l_2}\cdots 1^{k_t}2^{\infty}.$ Then
\begin{equation}\label{eq2.0}
\vert (S_{n} \varphi)(x)-(S_{n} \varphi)(x') \vert = \tau\vert \log
\frac{ \vert\pi^{-1}[12^{\beta-(n-k)}1]\vert}
{\vert\pi^{-1}[12^{\beta}1]\vert}\vert \leq \tau \log K.
\end{equation}

\em{Case 2.} \em  Consider $x$ and $x'$ such that
$\pi(x_i)=\pi(x'_i)$ for all $0\leq i\leq n-1$ and
$\pi(x_{n-1})=\pi(x'_{n-1})=1, $ or $x$ and $x'$ such that
$\pi(x_i)=\pi(x'_i)$ for all $i\geq 0.$  Then
\begin{equation}\label{eq3.0}
\vert (S_{n} \varphi)(x)-(S_{n} \varphi)(x')\vert=0,
\end{equation}
($\varphi(x_{n-1}\cdots)-\varphi(x'_{n-1}\cdots)=0$ for
$\pi(x_{n-1})=\pi(x'_{n-1})=1.$)

By (\ref{eq1.1}), (\ref{eq2.0}), and (\ref{eq3.0}), $\tau G\circ \pi
\in Bow (X, \sigma).$ Now use Lemma \ref{BandH}.
\end{proof}
\subsection{Proof of (2)(c).}\label{mixing}  Note that $G$ is a grid function but $G\circ \pi$ is not
necessarily a grid function. Let $\overline M_0= \pi^{-1}[1]=[1]
\text{ and } \overline M_n= \pi^{-1}[2^{n}1] \text { for } n\geq 1$.
Let $\overline Z=\pi^{-1}\{2^{\infty}\}.$ Recall Condition [C]. Let
$W$ be the two-sided shift of finite type which is the natural
extension of $\overline Z.$ Then we have $W =\Sigma_{T},$ and so
$\rho(W)=W\vert _{+}=\Sigma_{T}\vert_{+}=\overline Z$ by Condition
[C](1). Thus $\p=\{\rho(W), \overline M_0, \overline M_1, \overline
M_2,\cdots \}$ is a partition of $X$ for which $G\circ \pi \in \g(W,
\p).$

We show first that in this case for any $0\leq \tau < 1,$ $\tau
G\circ \pi$ satisfies the RPF condition. Let $\tau=\alpha/({\alpha
+1}), \alpha >0.$ Define $\varphi= (\alpha/({\alpha +1}))G\circ
\pi.$ Since $\varphi \in \g(W,\p)$, it is enough to show that
$P_{X}(\varphi)>h_\text{top}(W)$. Note that $h_\text{top}(W)=h_\text{top}(\overline Z)=0$ because
$h_\text{top}(X_B)=0.$ By Formulas (\ref{key1}) and (\ref{key2}), if we
denote the Shannon-Parry measure by $\mu_{\text {max}}$, we have
that
\begin{equation*}
P_{X}(\varphi)\geq \frac{1}{\alpha
+1}h_{\mu_{\text{max}}}(\sigma_X)=\frac{1}{\alpha+1}h_\text{top}(X)>0.
\end{equation*}
Hence $\varphi$ satisfies the RPF condition by Theorem \ref{MPkey}.
Let $\nu, h, \lambda$ be given by the RPF condition. By \cite{H1},
the unique equilibrium state $\mu$ for $\varphi$ is $\mu=\nu h$. If
we let
$$\overline \varphi=({\alpha}/(\alpha +1)) G\circ \pi +\log h-\log
h\circ \sigma -\log \lambda,$$ then $\mu$ is the $g$-measure for
$g=e^{\overline \varphi}$ \cite{H1}. Since $ L^{n}_{\overline
\varphi}f \rightarrow \mu(f)$ uniformly for all $f\in C(X)$
\cite{H1}, by using the same arguments as in the proof of Theorem
3.2 \cite{Wg}, we conclude that $(\sigma, \mu)$ is an exact
endomorphism, hence strongly mixing.

\section {Main Result-Part 2}\label{Main2}
We next consider the case when $h_\text{top}(X_{B})>0,$ under Setting
(A).
\begin{theorem}\label{thm2}
Let Setting (A) in Section \ref{mainthms} hold. Suppose that
$h_\text{top}(X_{B})>0$ and that the following two conditions hold:\\
(i)$'$ There exists $a> 1$ such that
\begin{equation*}\label{limquotient2}
\lim_{n \to\infty} \frac{\vert \pi^{-1}[12^{n-1}1]\vert}{\vert
\pi^{-1}[12^{n}1]\vert}=\frac {1}{a}.
\end{equation*}
(ii)$'$
\begin{equation*}\label{-topofXB2}
h_\text{top}(X_B)=\limsup_{n \to\infty}\frac{\log \vert
\pi^{-1}[2^{n}]\vert}{n}=\log{a}.
\end{equation*}
Then \\

\renewcommand{\theenumi}{\alph{enumi}}
\begin{enumerate}
\item There exists a compensation function $G\circ \pi \in C(X)$
such that $G\in C(Y).$
\item $\tau G\circ \pi$ has a unique equilibrium state $\mu$.
\item Under Condition [C] (see page \pageref{conditionC}), $(\sigma, \mu)$ is exact, hence strongly mixing.
\item The unique equilibrium state $\mu$ is Gibbs if and only if there exist $K_1, K_2>0$
such that $$ K_1\leq \frac{a^n}{\vert \pi^{-1}[12^{n}1]\vert}\leq
K_2$$ for all  $n \geq 1.$
\item If the unique equilibrium state $\mu$ is Gibbs, then the natural extension of
$(\sigma,\mu)$ is isomorphic to a Bernoulli system.

In particular, if $B$ is irreducible with $h_\text{top}(X_B)>0,$ we have
the
following: If (i)$'$ is satisfied, then (a),(b),(c) above hold and

\item the unique equilibrium state $\mu$ is
Gibbs
and the natural extension of $(\sigma,\mu)$ is isomorphic to a Bernoulli system.

\end{enumerate}
\end{theorem}

\begin{remarks}\label{remarkmain}
\textnormal{
1. The hypothesis that $X$ be topologically mixing is not necessary
for (a).\\
2. If $B$ is primitive, (i)$'$ is automatically satisfied.} \\
\end{remarks}
\begin{proof}
We give only the outlines of the proofs because they are similar to
those of Theorem \ref{thm1}. Note first that we no longer have
$G(2^{\infty})=0$ because $h_\text{top}(X_B)>0$. Define $G: Y\rightarrow
\R$ as in Formula (\ref{comfun}), but replacing 0 at $2^{\infty}$ by
$-\log a.$ We modify \em Case 1-2 \em of \ref{maincase1-1} and \em
Case 1-4 \em of \ref{maincase1-3} to take into account that $a>1$.
We find upper bounds for ${\vert \pi^{-1}[12^{t}1]\vert \vert
\pi^{-1}[12^{k_l}]\vert}/{\vert \pi^{-1}[12^{k_{l} +t}1]\vert}$ by
using (i)$'$ and (ii)$'$. For (i)$'$, fix $\epsilon>0$ small enough
so that $\epsilon +{1}/{a}<1.$ Then there exists $N \in \N$ such
that
$$ (\frac{1}{a}-\epsilon)\vert\pi^{-1}[12^{t+1}1]\vert \leq
\vert\pi^{-1}[12^{t}1]\vert\leq
(\frac{1}{a}+\epsilon)\vert\pi^{-1}[12^{t+1}1]\vert $$ for all
$t\geq N.$ Therefore, in the case when $t=t(y,n)\geq N$ and
$k_l=k_{l}(y,n)\geq 1,$ we have that
\begin{equation*} \label{ratioo}
\frac
{\vert\pi^{-1}[12^{t}1]\vert}{\vert\pi^{-1}[12^{t+k_l}1]\vert}\leq
(\frac{1}{a}+\epsilon)^{k_{l}}.
\end{equation*}

For (b), by Corollary \ref{unique}, it is enough to prove that
$(-1/(1+\alpha))(G+\log a)$ has a unique equilibrium state. Since
$(-1/(1+\alpha))(G +\log a) \in \g(Z,\p),$ where $Z, \p$ are defined
as in Section \ref{pfu0}, it is enough to show that
\begin{equation}\label{aun} P_{Y}(-\frac{1}{\alpha +1}(G+\log a))>h_\text{top}(Z)(=0).
\end{equation}
Using (\ref{key1}), (\ref{key2}), and $h_\text{top}(X)>h_\text{top}(X_B),$ we
get (\ref{aun}).

For (c), we show $\tau(G+\log a)\circ\pi$ satisfies the
RPF condition for $0\leq \tau<1,$ as in Section \ref{mixing}.

For (d), we replace $G(2^{\infty})=0$ by $G(2^{\infty})=-\log a$ and
repeat the arguments as in the proof of (2)(d) in Section \ref{WH}.

Conclusion (e) follows by Lemma \ref{BandH}.

Next consider the case when $B$ is irreducible. Then (ii)$'$ is
satisfied. (i)$'$ is satisfied if $B$ is primitive. To see this, let
$$V_1=\{c \in\{2, 3,..., r\}: A_{1c}\neq 0\} \text{ and } V_2=\{d \in\{2,
3,..., r\}: A_{d1}\neq 0\}.$$ Since ${\vert
\pi^{-1}[12^{n}1]\vert}=\sum_{c\in V_1, d\in V_2}(B_{c, d})^{n-1},$
apply the Perron-Frobenius Theorem to find the growth rate of
$(B_{c, d})^{n-1}.$ For (f), we use the following consequence of the
Perron-Frobenius Theorem.
\begin{lemma}
Let $p$ be the period of the irreducible matrix $B.$  There exist
$N, A_1, A_2>0$ such that
$$ A_1 a^{np+l}\leq \vert \pi^{-1}[12^{np+l+1}1]\vert \leq A_2
a^{np+l}$$ for all $l\in \{0,1,\cdots p-1\}, n\geq N.$
\end{lemma}
\noindent Now we apply Lemma \ref{claim1} and  Lemma \ref{BandH}.
\end{proof}

\section{Examples for Theorem \ref{thm1} and Theorem \ref{thm2}}
\label{fexamples} 
\begin{figure}[!h]
\centering
\scalebox{.4}{\includegraphics{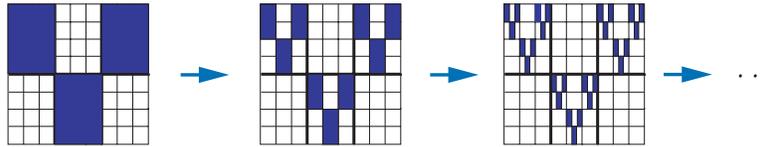}}\caption{NC
carpet $K(T,R)$ in Example \ref{ap2}}\label{NC1}
\end{figure}

\begin{figure}[!h]
\centering
\scalebox{.4}{\includegraphics{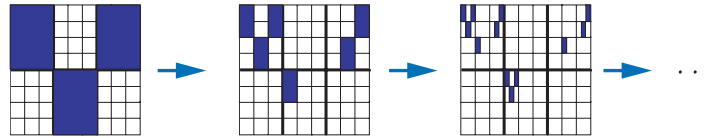}}
\caption{SFT-NC carpet $K(T,R,A)$ in Example \ref{ap2}}\label{SFTE1}
\end{figure}

We will apply Theorem \ref{thm1} and Theorem
\ref{thm2} to study SFT-NC carpets.
\begin{example} \label{ap2}
\end{example}
Let $T$ be the toral endomorphism given by $T(x,y)=(3x \textnormal{
mod }1, 2y \textnormal{ mod }1).$ Let $$\p=\{ [\frac{i}{3}, \frac
{i+1}{3}]\times[\frac{j}{2}, \frac{j+1}{2}]: 0\leq i \leq 2, 0\leq
j\leq 1\}$$ be the natural Markov partition for $T.$ Let $R=\{(1,0),
(0,1),(2,1)\}.$ Then we get the NC carpet $K(T, R)$ given in Figure
\ref{NC1}. Let $A$ be the transition matrix among the members of $R$
given by
\begin{equation*}
A=\left [\begin{array}{rrr}
                  0 &1& 0 \\
                  1 & 1&1 \\
                  1&0&1\\
                  \end{array}\right].
\end{equation*}
Then the SFT-NC carpet $K(T,R,A)$ is defined by
\begin{equation*}
K(T,R,A)=\{(\sum_{k=1}^{\infty}\frac{x_k}{3^{k}},
\sum_{k=1}^{\infty} \frac{y_k}{2^{k}}): (x_k, y_k)\in R, A_{(x_k,y_k)(x_{k+1}, y_{k+1})}=1 \textnormal{ for all } k \}.
\end{equation*}
(See Figure \ref{SFTE1}.)


Let $X=\Sigma^{+}_A$ and let $\pi$ be the factor map $\pi(x_k,
y_k)=y_k.$ Set $Y=\pi(X)$. Then $Y$ is a shift of finite type on the
two symbols $0,1$ with the transition matrix
\begin{equation*}
C=\left [\begin{array}{rr}
                  0 &1  \\
                  1 & 1 \\
                  \end{array}\right].
\end{equation*}
If we denote $(1,0)$ by 1, $(0,1)$ by 2, and $(2,1)$ by 3, then
$\pi(1)=0,\pi(2)=\pi(3)=1.$  Let $\alpha = \log_{2}3-1.$ Then the
transition matrix $B$ of symbols in $\pi^{-1}\{1\}$ is given by
\begin{equation*}
B=\left [\begin{array}{rr}
                  1 &1  \\
                  0 & 1 \\
                  \end{array}\right].
\end{equation*}

\bigskip
\begin{center}
\xymatrix{
& & & &     2 \ar@(ur,dr)[] \ar[dd]  \\
& & 1\ar@{<->}[urr] \ar@{<-}[drr] & & & & \ar [r]^\pi & & 0\ar@{<->}
[r] & 1 \ar@(ur,dr) \\
& & & &    3 \ar@(ur,dr)[]            }
\end{center}
\bigskip
\begin{figure}[h]
{} \caption{$X,Y,$ and $\pi$ in Example \ref{ap2}} \label{fig1}
\end{figure}

Now we apply Theorem \ref{thm1}. Since we have $\vert
\pi^{-1}[01^{n}0]\vert=n$ and $h_\text{top}(X_B)=0,$ (2)(i) and (ii) are
satisfied. There is a saturated compensation function $\com, $ where
$G: Y\rightarrow \R$ is defined by

\begin{equation*}
G(y)= \begin{cases} \log ((n-1)/n) & \textrm{
 if }y\in[1^{n}0], n\geq 2 \\
0 & \textrm{
 if }y\in [0] \cup [10] \cup \{1^{\infty}\}.
\end{cases}
\end{equation*}

$\acom$ has a unique equilibrium state, $\overline
\mu.$ Therefore, the corresponding measure $\mu$ on $K(T,R,A)$ is
the unique measure of full Hausdorff dimension. Since Condition [C]
is satisfied, $(T,\mu)$ is exact, and hence strongly mixing. The
unique measure $\mu$ is not Gibbs, because $\sup_{n} \vert
\pi^{-1}[01^{n}0]\vert=\infty.$ The RPF condition is satisfied, and
by Corollary \ref{HDK},
$$\dim_{H}K(T,R,A)=\frac{P(\acom)}{\log 2}=\log_{2} \lambda,$$
where $\lambda$ is the spectral radius of $L_{\acom}$.

We now approximate $P(\acom),$ which gives an approximation of
$\lambda.$  Let $\beta =\alpha /(\alpha +1)$ and $\varphi=\beta
G\circ\pi.$ Consider a lower bound for
$$\sum_{x_0
x_1 \cdots x_{3m-1}} \exp(\sup_{y  \in [x_0 x_1 \cdots x_{3m-1}]}
{(S_{3m} \varphi)(y)}).$$

First note that the number of allowable cylinder sets of the form
$[x_0 x_1 \cdots x_{m-1}]$ in $\Sigma^{+}_A$ of length $m$ is
$3\cdot 2^{m-1}.$  We consider the cylinder set $[x_0 x_1\cdots
x_{3m-1}]$, where 221, 331 or 231 appears $k$ times in total
somewhere in $x_0 x_1 \cdots x_{3m-1}.$ So there are $i_1, i_2,
\dots, i_k$ such that $x_{i_l}x_{i_{l}+1}x_{i_{l}+2}\in \{221, 331,
231\}$ for $l=1,2,\dots,k, 1\leq k\leq m.$ If $y \in [x_0 x_1\cdots
x_{3m-1}]$ and if 221, 331, or 231 appears in $x_0 x_1\cdots
x_{3m-1}$ a total of $k$ times, then direct computation shows that
$$ (S_{3m} \varphi)(y)\geq \beta [k \log(1/2)+ ( 3m-3k)\log(2/3)].$$
 Using this, we have a lower bound
\begin{displaymath}
\exp(\sup_{y \in [x_0 x_1 \cdots x_{3m-1}]} {(S_{3m} \varphi)(y)})
\geq (\frac{2}{3})^{3m\beta}.
\end{displaymath}
This implies by direct computation that
\begin{displaymath}
P(\varphi)\geq\lim_{m \rightarrow \infty} \frac{\log [3\cdot
2^{3m-1}\cdot(\frac{2}{3})^{3m\beta}]}{3m}\geq \log
((\frac{2}{3})^{\beta}\cdot 2)
>0.
\end{displaymath}
Therefore, we get $\lambda\geq (2/3)^{\beta}\cdot 2.$
\label{cqapplication}

Next we show that the natural extension of $(T, \mu)$ is isomorphic
to a Bernoulli system. Since the unique measure is a $g$-measure, we
use Theorem \ref{CQg} with the approximation of $\lambda$ from
above.

Let $h \in C(X)$ be an eigenfunction of the Perron-Frobenius
operator with corresponding eigenvalue $\lambda.$ Then the unique
measure $\mu$ is a $g$-measure for $g=e^{\varphi +\log h -\log h
\circ \sigma-\log \lambda}.$ Note that $g$ satisfies the hypothesis
of Theorem \ref{CQg}.

Applying Theorem \ref{CQg}, it is enough to show that
\begin{displaymath}
\sum_{n=r}^{\infty} \prod_{i=r}^{n} (1-\frac{3}{2} \var_{i}(g))=
\infty \textnormal{ for some r} \geq 1,
\end{displaymath}
where
$$g(x)=\frac {e^{\varphi(x)}h(x)}{(h \circ \sigma)(x)
\lambda}.$$ We show this by the following three
steps. \\
{\em Step 1}:   Show that there is a $K_0 \geq 1$ such that
$\var_{i}(\log g)\leq \log ({i}/(i-1))^{2\beta}$ for $i\geq
K_0 +1$. \\
{\em Step 2}: Show that $\var_i(g)\leq (({i}/(i-1))^{2\beta}
-1)/\lambda$ for
$i\geq K_0 +1.$\\
{\em Step 3}: Show that there is $r\geq 1$ such that
 \begin{displaymath}
\sum_{n=r}^{\infty} \prod_{i=r}^{n} (1-\frac{3}{2} \var_{i}(g))=
\infty \textnormal{ for some r} \geq 1.
\end{displaymath}
{\em Step 3} follows from {\em Step 2} easily:\\
Suppose $\var_i(g)\leq (({i}/(i-1))^{2\beta}-1)/{\lambda}$ for
$i\geq K_0 +1.$ Then we have that
$$
1-\frac{3}{2} \var_k(g)\geq
1-\frac{3}{2}((\frac{k}{k-1})^{2\beta}-1) \frac{1}{\lambda}$$
 for $k\geq K_0 +1.$ Since $\lambda >{3}/{2}$ (by the
 approximation of $\lambda$ above), there exists $K^{'}_0\in \N$ such that
$$
 1-\frac{3}{2}((\frac{k}{k-1})^{2\beta}-1) \frac{1}{\lambda} \geq (\frac{k-1}{k})^{2\beta}$$ for
all  $k\geq K^{'}_0.$
 Now take $\overline {K_0}= \max \{K_0 +1,K^{'}_0\}$.  Then
$$\prod_{i=\overline K_0}^{n} (1-\frac{3}{2} \var_i(g))\geq (\frac{\overline
K_0-1}{n})^{2\beta}.$$ Therefore, since $2\beta<1,$
$$\sum_{n=\overline K_0}^{\infty} \prod_{i=\overline K_0}^{n}(1-\frac{3}{2} \var_{i}(g)
)\geq \sum_{n=\overline K_0}^{\infty}(\frac{\overline
K_0-1}{n})^{2\beta} =\infty.$$

{\em Step 2} also follows from {\em Step 1}.  Thus we show
{\em Step 1}. Let $\overline \varphi=\log g.$
\begin{equation*}
 \var_i (\overline\varphi) = \sup \{ \vert\varphi(x)-\varphi(y) + \log {\frac{h(x)h(\sigma y)}{h(\sigma
x)h(y)}}\vert: x_k=y_k, \textnormal{ for all } 0\leq k\leq {i-1}\}.
\end{equation*}
 By Theorem \ref{MPkey}, we know that $h$ is constant on
any cylinder set $B$ of length larger than or equal to $K_0$ such that $B\cap
\{2^k 3^{\infty}, 2^{\infty} : k\geq 0\} =\emptyset $ for some
$K_0.$ Define $M_{\infty}=\{2^{k}3^{\infty}, 2^{\infty} :k\geq 0\}.$
Use this fact to show that $\var_i( \overline \varphi)\leq \log
({i}/(i-1))^{2\beta}$ for $i=K_0 +1.$ Using
the properties of $h$, we have the following three cases. \\
\renewcommand{\theenumi}{\Roman{enumi}}
\renewcommand{\theenumii}{\alph{enumii}}
\renewcommand{\labelenumii}{\theenumii.}
\begin{enumerate}
    \item $[x_0 x_1 \cdots x_{K_0}] \cap
           M_{\infty}=\emptyset$ and
           $[x_1 x_2 \cdots x_{K_0}] \cap M_{\infty}=\emptyset$\\

    \item $[x_0 x_1 \cdots x_{K_0 }] \cap M_{\infty} \ne
           \emptyset$ and
           $[x_1 x_2 \cdots x_{K_0}] \cap M_{\infty} \ne \emptyset$\\
            In this case, $[x_0 x_1\cdots x_{K_0}]\in [2^{K_{0} +1}] \cup
            [3^{K_{0}+1}]\cup [2^{k_1}3^{k_2}],\quad k_1 +k_2=K_0
            +1.$\\

    \item $[x_0 x_1 \cdots x_{K_0}] \cap M_{\infty}
          =\emptyset$ and
          $[x_1 x_2 \cdots x_{K_0}] \cap M_{\infty} \ne \emptyset$\\
           In this case, $[x_0 x_1\cdots x_{K_0}]\in [12^{K_0}] \cup
           [12^{k_1}3^{k_2}]\quad \textnormal{with } k_1 +k_2=K_0, k_1\geq 1.$ \\
\end{enumerate}

Clearly for cylinder sets of type (I), we have that  $h(x)=h(y)$ and
$h(\sigma x)=h(\sigma y).$ Therefore,

\begin{eqnarray*}
\var_{K_0 +1}(\overline
\varphi)\vert_{\textnormal{cylinder sets of type I}} & \leq & \sup\{\vert
\varphi (x)-\varphi(y) \vert: x_k=y_k \textnormal{ for all }0 \leq
k\leq K_0\}\\ & \leq & \log (\frac{K_0 +1}{K_0})^{ 2\beta}.
\end{eqnarray*}

For cylinder sets of type (II), we can similarly show that
$$\var_{K_0 +1}(\overline \varphi)\vert_{\textnormal{cylinder sets
of type II}}\leq\log (\frac{K_0 +1}{K_0})^{2\beta},$$ by showing that
\begin{equation}\label{ineqc}
\vert \varphi(x)-\varphi(y) + \log {\frac{h(x)h(\sigma y)}{h(\sigma
x)h(y)}}\vert  \leq \  \log (\frac{K_0 +1}{K_o})^{2\beta}
\end{equation}
 for $x,y\in [x_0 x_1\cdots x_{K_0}]\in [2^{K_0+1}] \cup
[3^{K_0+1}]\cup [2^{k_1}3^{k_2}].$

 We will need to consider the following
cylinder sets of type (II).
\renewcommand{\theenumi}{\roman{enumi}}
\renewcommand{\theenumii}{\alph{enumii}}
\renewcommand{\labelenumii}{\theenumii.}
\begin{enumerate}
\item $ x,y\in [2^{K_0+1}],$
     \begin{enumerate}
     \item $x\in [2^n 1], y\in [2^m 1], \quad  n,m\geq K_{0} +1;$
     \item $x\in [2^n 1], y=2^{\infty}, \quad  n\geq K_{0} +1 ;$
     \item $x\in [2^n 1], y\in [2^m 3^{k}1], \quad  n,m\geq K_{0} +1, k \geq 1;$
     \item $x\in [2^n 1], y=2^m 3^{\infty}, \quad  n,m\geq K_{0} +1;$
     \item $x=2^{\infty}, y=2^n 3^{\infty},  \quad n\geq K_{0} +1;$
     \item $x=2^{n}3^{\infty}, y=2^m 3^{\infty},  \quad n, m\geq K_{0} +1;$
     \item $x\in [2^{n} 3^{k}1], y= 2^{\infty}, \quad n\geq K_{0} +1, k
     \geq 1 ;$
     \item $x\in [2^{n} 3^{k}1], y= 2^{m}3^{\infty}, \quad n,m\geq K_{0} +1, k \geq 1 ;$
     \item $x\in [2^{n} 3^{k_1}1], y\in [2^{m} 3^{k_2}1], \quad n, m\geq K_{0} +1, k_1, k_2\geq 1 ;$
     \end{enumerate}
\item $ x,y\in [3^{K_0+1}];$
     \begin{enumerate}
     \item $x\in [3^n 1], y\in [3^m 1], \quad n,m\geq K_{0} +1;$
     \item $x\in [3^n 1], y=3^{\infty} \quad n \geq K_{0} +1; $
     \end{enumerate}
\item $ x,y\in [2^{k_1} 3^{k_2}] \quad k_1+k_2=K_0+1, k_1, k_2\geq 1;$
      \begin{enumerate}
     \item $x\in [2^{k_1} 3^{l}1], y\in
          [2^{k_1} 3^{m}1], \quad l, m  \geq k_2 ;$
     \item $x\in [2^{k_1} 3^{l}1], y=2^{k_1} 3^{\infty}, \quad l \geq k_2.$
\end{enumerate}
\end{enumerate}
The idea to show (\ref{ineqc}) is to write $h(x), (h\circ \sigma) (x),
h(y), (h\circ \sigma) (y)$ as infinite series by using the property
that $h$ is an eigenfunction of the Perron-Frobenius operator with
corresponding eigenvalue $\lambda$.  For simplicity, we only
consider the case (i)b with $n=K_{0}+1$. Let $x\in [2^{K_0+1} 1], y=2^{\infty}.$ We
show that
$$
\vert \overline \varphi(x)- \overline \varphi(y) \vert \leq \log
(\frac{K_0 +1}{K_0})^{2\beta}.$$ Similar arguments work for the
other cases. Since
$$\vert \overline
\varphi(x)-\overline \varphi(y) \vert=\vert \log(\frac{K_0 }{K_0
+1})^{\beta}+ \log \frac {h(x)}{h(\sigma x)}\vert,$$ it is enough to
show that
$$1 \leq \frac{h([2^{K_0
+1}1])}{h([2^{K_0}1])}\leq \log(\frac{K_0 +1 }{K_0})^{\beta}.$$ Using
the Perron-Frobenius operator, we have that
 $$h([2^{K_0}1])=
\frac{1}{\lambda}[h([12^{K_0}1])+h([2^{K_0 +1}1])(\frac{K_0 }{K_0
+1})^{\beta}].$$
 Similarly, $$h([2^{K_0 +1}1])=
\frac{1}{\lambda}[h([12^{K_0 +1} 1])+h([2^{K_0 +2}1])(\frac{K_0 +1
}{K_0 +2})^{\beta}].$$ Note that we have $h([12^{K_0} 1])=h([12^{K_0
+1} 1]),$ because $[12^{K_0}] \cap M_{\infty}=\emptyset.$ Let
$C=h([12^{K_0} 1])$.  Repeating this argument, we have that
$$h([2^{K_0}1])=\frac {C}{\lambda}
(\sum_{i=0}^{n-1} \frac{{K_0}^{\beta}}{{\lambda}^{i}({K_0}
+i)^{\beta}}) + \frac{1}{{\lambda} ^n}(\frac {K_0}{{K_0}
+n})^{\beta} h([2^{K_0 +n}1])\textnormal { for all } n\geq 1.$$
Letting $n \rightarrow\infty$, we have that$$h([2^{K_0}1])=\frac
{C}{\lambda}\sum_{i=0}^{\infty}\frac{{K_0}^{\beta}}{{\lambda^{i}}{({K_0}
+i)}^{\beta}}\quad.$$ Similarly, we have $$h([2^{K_0 +1}1])=\frac
{C}{\lambda}\sum_{i=0}^{\infty}\frac{{(K_0
+1)}^{\beta}}{{\lambda^{i}}{(K_0 +1+i)}^{\beta}}\quad.$$ This
implies that
$$\frac{h(x)}{h(\sigma x)} = \frac{h([2^{K_0
+1}1])}{h([2^{K_0}1])}\leq \log(\frac{K_0 +1 }{K_0})^{\beta}.$$ It is
easy to see by direct computation, using the infinite series, that
 $$1\leq \frac{h(x)}{h(\sigma x)}.$$ Now the
result follows.

For the cylinder sets of type (III), we use similar arguments; find
the value of $h$ on a cylinder set by using the Perron-Frobenius
operator.\\

\noindent From now on we give symbolic dynamical examples without
presenting the particular SFT-NC carpets from which they arise. Note
that the carpets always exist by defining $T$ appropriately.

\begin{figure}[h]
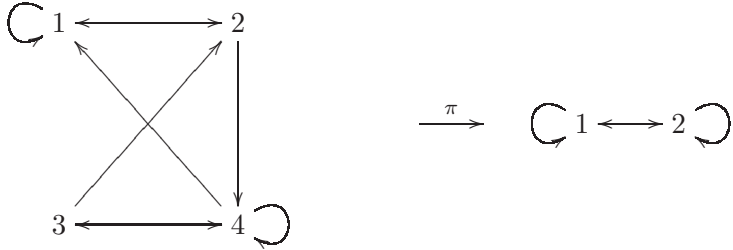

{} \caption{$X,Y,$ and $\pi$ in Example \ref{ap1}} \label{fig11}
\end{figure}

\begin{example}\label{ap1}
\end{example}
Let $X\subset \{1, 2, 3, 4\}^{\N}$ and $Y= \{1, 2\}^{\N}$ be the
shifts of finite type determined by the transitions given by Figure
\ref{fig11}. Define $\pi$ by $\pi(1)=1,$ $\pi(2)=\pi(3)=\pi(4)=2.$

\bigskip
\begin{center}
\xymatrix{
& & 1 \ar@(ul,dl) \ar@{<->}[rr] & &   2 \ar@{->}[dd] \ar@{<-}[ddll]\\
& & & & & & \ar [r]^\pi  & & 1 \ar@(ul,dl) \ar@{<->}
[r] & 2 \ar@(ur,dr)\\
& & 3  \ar@{<->}[rr]   & & 4 \ar@{->}[uull] \ar@(ur,dr) }
\end{center}
\bigskip

Then we have \\
\begin{displaymath}
A=\left [\begin{array}{rrrr}
                  1 &1& 0&0 \\
                  1 & 0&0 &1\\
                  0&1&0&1\\
                  1&0&1&1\\
                  \end{array}\right]
\textnormal{ and }
B=\left [\begin{array}{rrr}
                  0 &0&1  \\
                  1 & 0 &1\\
                  0&1&1\\
                  \end{array}\right].
\end{displaymath}
$A$ and $B$ are both primitive. Thus there is a compensation
function $G\circ \pi\in Bow(X,\sigma).$ Condition [C] is satisfied. Therefore, $(\sigma, \mu)$ is
exact, and hence strongly mixing.

\section{Some extensions of Theorem \ref{thm1} and Theorem
\ref{thm2}}\label{extensions} \indent We generalize Theorem
\ref{thm1} and Theorem \ref{thm2} to the case when the alphabet of
$Y$ can have more than two symbols. In this section, for an
allowable word  $y_{0}y_{1}\cdots y_{n-1}1$ of length $(n+1)$ in $Y$
with $y_{i}\neq 1$ for $0\leq i\leq n-1,$ we define $\vert
\pi^{-1}[1y_{0}y_{1}\cdots y_{n-1}1]\vert =\vert
\pi^{-1}[y_{0}y_{1}\cdots y_{n-1}1]\vert\geq 1 $ if
$\pi^{-1}[1y_{0}y_{1}\cdots y_{n-1}1] ={\emptyset}.$\\

\noindent \textbf{\underline{Setting (B)}}\\
Fix $k=3,4,\cdots.$ Let $r_2,\cdots,r_k\in \N$, let $a^{j}_{i}
  (j=2,\cdots, k; i=1,\cdots, r_k)$ be symbols, and let $X \subset
\{1,a^2_{1},\cdots,a^{2}_{r_2},a^3_{1},\cdots,a^{3}_{r_3},
\cdots,a^k_{1},\cdots, a^{k}_{r_k}\}^{\N}$ be a topologically mixing
shift of finite type with positive entropy. Let $Y= \{1,2,\cdots,
k\}^{\N}$ or $Y\subset \{1,2,\cdots,k\}^{{\N}}$ be a shift of finite
type with positive entropy and $\pi : X \rightarrow Y$ a one-block
factor map such that $\pi^{-1} \{1\}=\{1\}$ and $\pi^{-1}
\{i\}=\{a^i_{1},\cdots,a^{i}_{r_i}\}$ for $2\leq i\leq k$. Let $A$
be the transition matrix of $X$, let $B$ be the submatrix of $A$
corresponding to the indices $a^2_{1},\cdots, a^{k}_{r_k}$ (giving
the transitions among symbols in $\pi^{-1}\{2,\cdots,k\}$) and let
$B'$ be the transition matrix of $\{2,3,\cdots,k\}.$ Assume that for
all $n$, $B^n\neq 0.$ Denote by $X_B$ the shift of finite type
determined by $B$ and by $Y_{B'}$ the shift of finite type
determined by $B'$. Let $0\leq \tau<1.$\\

\noindent \textbf{Condition [C$'$]}\\
(1) If $\Sigma_{B'}$ is the two-sided shift of finite type on
$\{2,3,\cdots,k\}$ determined by $B'$, and $\Sigma_{B'}\vert _{+}$
is the projection of $\Sigma_{B'}$ onto a one-sided shift of
finite type, then $\Sigma_{B'}\vert _{+}=Y_{B'}.$\\
(2) The ratio
$$\frac{\vert \pi^{-1}[1y_{1}y_{2}\cdots
y_{n-1}1]\vert}{\vert \pi^{-1}[1y_{0}y_{1}\cdots y_{n-1}1]\vert}$$
is constant in $n$, for all allowable words $y_{0}\cdots y_{n-1}1$
of length $(n+1)$ in $Y$ such that no $y_i= 1.$\\
(3) Let $S=\{s\in \{2,\cdots, k\}:  s1 \text{ is allowed in } Y\}$.
For any $a \in \{2,\cdots,k\}$, there exist $s\in S, K\in \N, y_0 \cdots y_{K-1},$ a word of length $K$ in
$Y_{B'}$ (where $K$ is independent of $a$ and $s$)
such that $ay_0\cdots y_{K-1}s$ is allowable in $Y_{B'}.$
\begin{theorem}\label{prop1}
Let Setting (B) hold and suppose the following hypothesis I is
satisfied.

I.  $h_\text{top}(X_{B})=0$ and the following condition holds: For each $y\in
Y_{B'},$
\begin{equation*}\label{gcon1}
\lim \frac{\vert \pi^{-1}[1y_{1}y_{2}\cdots y_{n-1}1]\vert}{\vert
\pi^{-1}[1y_{0}y_{1}\cdots y_{n-1}1]\vert}=1,
\end{equation*}
where the limit is taken along the sequence of $n\rightarrow \infty$
for which $y_{n-1}1$ is allowable in $Y.$ Then
\renewcommand{\theenumi}{\alph{enumi}}
\renewcommand{\theenumii}{\alph{enumii}}
\renewcommand{\labelenumii}{\theenumii.}
     \begin{enumerate}
     \item  There exists a compensation function $G\circ \pi \in C(X)$
    such that $G \in C(Y)$.
     \item Under Condition [C$'$] above, $G$ is a grid function.
     \item If $G$ is a grid function constructed in (b), $\tau G\circ \pi$ has a unique
      equilibrium state.
     \end{enumerate}
\end{theorem}
\begin{proof}
For (a), define $G: Y\rightarrow \R$ as
\begin{equation*} G(y)= \begin{cases}
\log ({\vert \pi^{-1}[1y_{1}y_{2}\cdots y_{n-1}1]\vert}/{\vert
\pi^{-1}[1y_{0}y_{1}\cdots y_{n-1}1]\vert}) & \textrm{
if }y\in[y_{0}y_{1}\cdots y_{n-1}1],\\
& y_i\neq 1, n\geq 2, \\
\log ({1}/{\vert \pi^{-1}[1y_{0}1]\vert}) & \textrm{
if }y\in[y_{0}1], \\
0 & \textrm{ if }y\in [1]\cup {Y_{B'}}.
\end{cases}
\end{equation*}
We can show that (\ref{Nice eq}) holds for all $y\in Y\setminus
E_1$.

For (b), note first that $G$ is not always a grid function when $Y$
has more than two symbols. Define $M_0=[1].$ For $n\geq 1,$ let $M_n= \cup[y_{0}\cdots y_{n-1}1]$, where the union is taken over all allowable words $y_{0}y_{1}\cdots y_{n-1}1$ of length $n+1$ in $Y$ such that $y_{i}\neq 1$ for $0\leq i\leq n-1.$ Let
$\p=\{\rho(\Sigma_{B'}), M_{0},M_1, M_2,\cdots\}.$ Then $G\in
\g(\Sigma_{B'}, \p).$

For (c), note that we cannot use Theorem \ref{G(2)}, and so we
modify the proof of 2(b) in Section \ref{pfmain12a}, using
$h_\text{top}(Y_{B'})=h_\text{top}(\Sigma^{+}_{B'})=h_\text{top}(\Sigma_{B'})=0,$
and $G\equiv0$ on $Y_{B'}.$
\end{proof}
\begin{theorem} \label{prop2}
Let Setting (B) hold and suppose the following hypothesis II is
satisfied.

II. $h_\text{top}(X_{B})=\log a, a>1,$ and the following condition holds: For
each $y\in Y_{B'},$
\begin{equation*}\label{gcon1}
\lim \frac{\vert \pi^{-1}[1y_{1}y_{2}\cdots y_{n-1}1]\vert}{\vert
\pi^{-1}[1y_{0}y_{1}\cdots y_{n-1}1]\vert}=\frac{1}{a} ,
\end{equation*}
where the limit is taken along the  sequence of $n \rightarrow
\infty$ for which $y_{n-1}1$ is allowed in $Y$.

Then
\renewcommand{\theenumi}{\alph{enumi}}
\renewcommand{\theenumii}{\alph{enumii}}
\renewcommand{\labelenumii}{\theenumii.}
    \begin{enumerate}
    \item  There exists a compensation function $G\circ \pi \in C(X)$
    such that $G \in C(Y)$.
    \item Under Condition [C$'$], $G+\log a$ is a grid function.
    \item If $G+ \log a$ is a grid function constructed in
    (b), $\tau G\circ \pi$ has a unique
    equilibrium state.
    \end{enumerate}
\end{theorem}

\begin{remark}
\textnormal{The hypothesis II seems to hold only in cases when $X_B$ is
reducible.  We study the case when $X_B$ is reducible in the next
two theorems, using different approaches to prove the existence of a
saturated compensation function and uniqueness for it.}
\end{remark}
\begin{proof}
For (a), define $G$ as in Theorem \ref{prop1} by replacing 0 by
$-\log a$ on $Y_{B'}.$
For each $n=1,2,\cdots$ and $y\in Y$, denote $F_{n}(y)$ by a set consisting of exactly one point from each nonempty
cylinder $[x_{0}\cdots x_{n-1}]\subset \pi^{-1}[y_{0}\cdots
y_{n-1}].$
We prove that for all $y\in Y\setminus
E_1,$ we have
\begin{equation}\label{upperb}
\limsup_{n\rightarrow\infty}\frac{1}{n}[\log(e^{(S_nG)(y)}\cdot \vert F_{n}(y)\vert)]\leq 0
\end{equation}
and
\begin{equation}\label{lowerb}
\limsup_{n\rightarrow\infty}\frac{1}{n}[\log(e^{(S_nG)(y)}\cdot \vert D_{n}(y)\vert)]\geq 0.
\end{equation}
For $y\in Y_{B'}\cup \{y\in Y\setminus Y_{B'} : y=y_{0}\cdots y_{p-1}b \textnormal{ for some } p\geq 1, b\in Y_{B'}\},$ we use a theorem of Petersen and Shin \cite{PS}. For such $y$, instead of showing (\ref{lowerb}), we show that
\begin{equation}\label{lowerb2}
\limsup_{n\rightarrow\infty}\frac{1}{n} \log[\Sigma_{x\in F_{n}(y)} e^{(S_n(G\circ \pi))(x)}]\geq 0
\end{equation}
by finding a lower bound of $\Sigma_{x\in F_{n}(y)} e^{(S_n(G\circ \pi))(x)}.$

For (b), Let $M_n$ and $\p$ be as in the proof of Theorem
\ref{prop1} (b). Then $G+\log a\in \g(\Sigma_{B'}, \p).$

For (c), we modify the proof of Theorem \ref{thm1}(2)(b). We show
first that if $-(1/(\alpha+1))G + (\alpha/(\alpha+1))\log
a$ has a unique equilibrium state, then $(\alpha/(\alpha+1))(G+\log
a)\circ \pi$ has a unique equilibrium state and hence so does $\tau
G\circ\pi,$ where $\tau=\alpha/(\alpha +1).$ Let $
\varphi=-(1/(\alpha+1))G + (\alpha/(\alpha+1))\log a $ and $\overline
\varphi=(\alpha/(\alpha+1))(G+\log a)\circ \pi.$ By Theorem \ref{CC},
$\mu$ is an equilibrium state of $\overline \varphi$ if and only if
$\pi\mu$ is an equilibrium state of $\varphi$ and $\mu$
is a relative equilibrium state of $G\circ \pi$ over $\pi\mu$. Let $\nu$ be the unique (ergodic) equilibrium state for
$\varphi$ and let $\mu$ be a preimage of $\nu$ with maximal relative
entropy.  Recall that $\mu$ is a preimage of maximal entropy if and
only if it is a relative equilibrium state of $G\circ \pi$ over
$\nu.$

Next we will show that $\nu([1])>0.$ If $\nu([1])=0,$ then
$\nu(Y_{B'})=1.$ Since $G(Y_{B'})=-\log a,$ using the variational
principle,
\begin{equation}
\begin{split}\label{P1}
P_{Y}(\varphi) &= h_{\nu}(\sigma_{Y})+ \int
_{Y_{B'}}(-\frac{1}{\alpha +1}G
+\frac{\alpha}{\alpha+1}\log a) d\nu \\
&=h_{\nu}(\sigma_{Y})+ \log a\geq \log a .
\end{split}
\end{equation}
Since $\mu$ is an equilibrium state of $\overline{\varphi}$ and
$\mu(X_B)=1,$
\begin{equation} \label{P2}
\begin{split}
P_{X}(\overline{\varphi}) &= h_{\mu}(\sigma_{X})+
\int_{X_{B}}(\frac{\alpha}{\alpha +1}G\circ \pi
+\frac{\alpha}{\alpha+1}\log a) d\mu \\
&=h_{\mu}(\sigma_{X})\leq
h_\text{top}(X_{B})=\log a.
\end{split}
\end{equation}
Since $G\circ \pi$ is a compensation function, we have $
P_{Y}(\varphi)= P_{X}(\overline{\varphi}),$ and so by using
(\ref{P1}) and (\ref{P2}), we get $P_{Y}(\varphi)=
P_{X}(\overline{\varphi})= \log a.$ Now by Formula (\ref{key1}) and
$h_\text{top}(X)>\log a$, taking the Shannon-Parry measure $\mu_{max}$ on
$X$, we get
\begin{equation*}
\begin{split}
P_{X}(\overline{\varphi}) &=P_{X}(\frac{\alpha}{\alpha +1}G\circ
\pi)+\frac{\alpha}{\alpha+1}\log a\\
&= \frac{1}{\alpha +1}\sup_{\mu \in M(X,\sigma_{X})} \{h_{\mu}(\sigma_{X})+\alpha h_{\pi\mu}(\sigma_{Y})\}+\frac{\alpha}{\alpha+1}\log a\\
& \geq \frac{1}{\alpha
+1}h_{\mu_{\max}}(\sigma_{X})+\frac{\alpha}{\alpha+1}\log a=
\frac{1}{\alpha +1}h_\text{top}(X)+\frac{\alpha}{\alpha+1}\log a
>\log a.
\end{split}
\end{equation*}
This is a contradiction. Therefore, $\nu([1])>0$. By using the same
arguments as in Section \ref{pfu0}, we conclude that if $\varphi$
has a unique equilibrium state, then $\overline{\varphi}$ has a
unique equilibrium state.

Thus it is enough to show that $\varphi$ has a unique equilibrium
state. If we let $\tilde {\varphi}= -(1/(\alpha +1))(G + \log a),$
then $P_{Y}(\tilde {\varphi})\geq h_\text{top}(Y_{B'})$ by the
variational principle. Since $\tilde \varphi$ is a grid function, if
the strict inequality occurs, we are done by Theorem \ref{MPkey}. If
not, then by Corollary 3.8 in \cite{MP}, exactly one of the
following happens: (i) The set of equilibrium states for
$\tilde{\varphi}$ is the set of measures of maximal entropy for
$\sigma_{Y_{B'}}.$ or (ii) The set of equilibrium states is the
closed convex hull of measures of maximal entropy for
$\sigma_{Y_{B'}}$ and the unique (ergodic) equilibrium state $\nu$
for  $\tilde{\varphi}$ such that $\nu(Y_{B'})=0.$ We show that (i)
does not happen. Assume that a measure $\bar \nu$ of maximal entropy
for $\sigma_{Y_{B'}}$ is an equilibrium state for $\tilde{\varphi}.$
Let $\bar \mu$ be a preimage of $\bar \nu$ with maximal relative
entropy.  Then $\bar \mu$ is a relative equilibrium state of $G\circ
\pi$ over $\bar \nu.$ Therefore by Theorem \ref{CC} $\bar \mu$ is an
equilibrium state for $\acom$, and so
\begin{equation*}
\begin{split}
P_{X}(\frac{\alpha}{\alpha+1}G\circ \pi)&= \frac{1}{\alpha
+1}\sup_{\mu\in
M(X,\sigma_{X})}\{h_{\mu}(\sigma_{X})+\alpha h_{\pi \mu}(\sigma_{Y})\}\\
 &=
\frac{1}{\alpha +1}h_{\bar \mu}(\sigma_{X})+\frac{\alpha}{\alpha +1}
h_{\bar \nu}(\sigma_{Y})\\
&=\frac{1}{\alpha +1}h_{\bar \mu}(\sigma_{X})+\frac{\alpha}{\alpha
+1} h_\text{top}(Y_{B'}).
\end{split}
\end{equation*}
Using the definition of compensation function and $G(Y_{B'})=-\log
a,$ we get
\begin{equation*}
\begin{split}
P_{X}(\frac{\alpha}{\alpha+1}G\circ
\pi)&=P_{Y}(-\frac{1}{\alpha+1}G) =h_{\bar
\nu}(\sigma_{Y_{B'}})+\frac{1}{\alpha +1}\log
a\\&=h_\text{top}({Y_{B'}})+\frac{1}{\alpha +1}\log a.
\end{split}
\end{equation*}
Therefore, $h_{\bar \mu}(\sigma_{X})-\log a= h_\text{top}(Y_{B'}).$ Since
$\bar \nu(Y_{B'})=1,$ we have $\bar \mu(X_{B})=1,$ and so using
$h_{\bar \mu}(\sigma_{X})\leq \log a,$ we get  $h_\text{top}(Y_{B'})=0.$
Thus $P_{Y}(\tilde \varphi)=h_\text{top}(Y_{B'})=0.$ However, taking the
Shannon-Parry measure $\mu_{max}$ on $X$ and using Formula
(\ref{key1}) and $h_\text{top}(X)>\log a,$ we get

\begin{equation*}
\begin{split}
P_{Y}(\tilde \varphi) &=P_{X}(\frac{\alpha}{\alpha +1}G \circ \pi )-\frac{1}{\alpha+1}\log a \\
&=\frac{1}{\alpha +1}\sup_{\mu\in M(X,\sigma_{X})}\{h_{\mu}(\sigma_{X})+\alpha
h_{\pi \mu}(\sigma_{Y})\}-\frac{1}{\alpha +1}\log a
\\
& \geq \frac{1}{\alpha +1}(h_{\mu_{ \max}}(\sigma_{X})-\log a) =
\frac{1}{\alpha +1}(h_\text{top}(X)-\log a)
>0.
\end{split}
\end{equation*}
This is a contradiction. Therefore, a measure of maximal entropy for
$\sigma_{Y_{B'}}$ is not an equilibrium state for $\tilde \varphi.$
Therefore, (ii) happens. Then there is a unique equilibrium state
$\nu$ for $\tilde{\varphi}$ such that $\nu(Y_{B'})=0$.
\end{proof}

The following proposition follows from the preceding proof of (c).
\begin{proposition}\label{hopetobeuseful}
Fix $k=2,3,\cdots.$  Let $X \subset\{1,2,\cdots,r\}^{\N}$ be a
topologically mixing shift of finite type with positive entropy, $Y=
\{1,2,\cdots, k\}^{\N}$ or $Y\subset \{1,2,\cdots,k\}^{{\N}}$ a
shift of finite type with positive entropy, and $\pi : X \rightarrow
Y$ a one-block factor map. Suppose there exist a partition $\p=\{A,
A_1, A_2\cdots \}$ of $Y,$ a two-sided subshift $(Z,\sigma)$ such
that $\rho (Z)=A$ (see page \pageref{grid}),  $a_{n}\in \R$ for
$n\in \N$, and a function $G: Y\rightarrow \R$ defined by
\begin{equation*}
G(y)= \begin{cases}
 a_n & \textrm{ if }y\in A_n,\\
-h_\text{top}(\pi^{-1}(A)) & \textrm{ if }y\in A, \\
\end{cases}
\end{equation*}
such that $G+h_\text{top}(\pi^{-1}(A)) \in \g(\p, Z)$ and $G\circ\pi$ is
a saturated compensation function. Then for each $\alpha>0,$
$-(1/(\alpha+1))G$ has a unique equilibrium state and the preimages
of the unique equilibrium state with maximal relative entropy are
equilibrium states for $\acom.$
\end{proposition}
\begin{remark}
\textnormal{The conditions when the hypotheses above hold need to be studied.
Saturated compensation functions in Theorem \ref{thm1} and Theorem
\ref{thm2} are examples for Proposition \ref{hopetobeuseful}.}
\end{remark}
In the next two theorems, we consider the case when $B$ is
reducible.
\begin{theorem}\label{prop3}
Let $X \subset
\{1,a^{2}_{1},\cdots,a^{2}_{r_2},a^3_{1},\cdots,a^{3}_{r_3}\}^{\N}$
be a topologically mixing shift of finite type with positive
entropy, $Y= \{1,2,3\}^{\N}$ or $Y\subset \{1,2,3\}^{{\N}}$ a shift
of finite type with positive entropy, and $\pi : X \rightarrow Y$ a
one-block factor map such that $\pi^{-1} \{1\}=\{1\}$ and $\pi^{-1}
\{i\}=\{a^i_{1},\cdots,a^{i}_{r_i}\}$ for $i=2,3$. Suppose that 23
and 32 are not allowable words in $Y$. Let $B_2$ be the submatrix of
$A$ corresponding to the indices $a^2_{1},\cdots,a^{2}_{r_2}$
(giving the transitions among $\pi^{-1} \{2\}$) and let $B_3$ be the
submatrix of $A$ corresponding to the indices
$a^3_{1},\cdots,a^{3}_{r_3}$  (giving the transitions among
$\pi^{-1}\{3\}$). Denote by $X_{B_2}$ the shift of finite type
determined by $B_2$ and by $X_{B_3}$ the shift of finite type
determined by $B_3$. Suppose for all $n,$ $(B_i)^n\neq 0$ for each
$i=2,3$. Let $B'$ be the transition matrix of $\pi^{-1}\{2,3\}$.

\renewcommand{\theenumiii}{\alph{enumiii}}
\renewcommand{\theenumii}{\arabic{enumii}}
\renewcommand{\theenumi}{\arabic{enumi}}
\renewcommand{\labelenumii}{\theenumii.}
\begin{enumerate}
\item  \label{con1prop3} Suppose that for each $i=2,3$ there exists $b_i \geq 1$ such that
\begin{equation*}\label{gcon1}
\lim_{n \to\infty} \frac{\vert \pi^{-1}[1i^{n-1}1]\vert}{\vert
\pi^{-1}[1i^{n}1]\vert}=\frac{1}{b_i} \textnormal{ and }
h_\text{top}(X_{B_i})=\log b_i.
\end{equation*}
Then there exists a compensation function $G\circ \pi \in C(X)$ such
that $G\in C(Y)$.

\item \label{con2prop3} Suppose in addition that
$b_i >1$ for $i=2,3$ in (1) and there exist $K_{1}(i), K_{2}(i)>0$
such that
\begin{equation*}
K_{1}(i)\leq \frac{b^n_i}{\vert \pi^{-1}[1i^{n}1]\vert}\leq K_{2}(i)
\textnormal{ for all } n.
\end{equation*}

Then $G\circ \pi \in Bow(X,\sigma).$ Hence $\tau G\circ \pi \in
Bow(X,\sigma)$ for all $\tau\in \R.$

\end{enumerate}
\end{theorem}
\begin{remarks} \textnormal{1. If $B_2$ and $B_3$ are both primitive, then (\ref{con1prop3})
and (\ref{con2prop3}) are automatically satisfied. If $B_2$ and
$B_3$ are irreducible, we still need (1) for $G\circ\pi$ to be in
$Bow(X,\sigma)$. These cases are also covered by Theorem
\ref{arigatou}.\\
2. The hypothesis of (2) is applicable to some cases not covered in
Theorem \ref{arigatou} when $X_{B}$ is reducible and two distinct
communicating classes for $X_{B}$ are mapped to the same symbol.
(see Example \ref{newe}). \\
3. We can generalize Theorem \ref{prop3} for the case when the
matrix of $\pi^{-1} \{2,3,\cdots k\}$ is a direct sum of
matrices of $\pi^{-1}\{i\}$ for $i=2,3,\cdots k.$}\\
\end{remarks}
\begin{proof}
For (1), define $G: Y\rightarrow \R$ as
\begin{equation*}G(y)= \begin{cases}
 \log ({\vert \pi^{-1}[1i^{n-1}1]\vert}/{\vert \pi^{-1}[1i^{n}1]\vert}) &
\textrm{
 if }y\in[i^{n}1], n\geq 2,i=2,3, \\
\log ({1}/{\vert \pi^{-1}[1i1]\vert}) & \textrm{
 if }y\in[i1],i=2,3, \\
-\log b_i & \textrm{
 if }y=i^{\infty}, i=2,3,\\
0 & \textrm{
 if }y\in [1].
\end{cases}
\end{equation*}
For (2), we modify Lemma \ref{claim1}, noting that $\pi(x_i)$ in the
proof of the lemma is 1, 2, or 3.
\end{proof}

We now consider another particular setting in which $B$ is
reducible.
\begin{theorem}\label{arigatou}
Let $X \subset
\{1,a^2_{1},\cdots,a^{2}_{r_2},a^3_{1},\cdots,
a^{3}_{r_3}\}^{\N}$ be a
topologically mixing shift of finite type with positive entropy, $Y=
\{1,2,3\}^{\N}$ or $Y\subset \{1,2,3\}^{{\N}}$ a
shift of finite type with positive entropy, and $\pi : X \rightarrow
Y$ a one-block factor map such that $\pi^{-1} \{1\}=\{1\}$ and
$\pi^{-1} \{i\}=\{a^i_{1},\cdots,a^{i}_{r_i}\}$ for $i=2,3$.
Let $B$ be the transition submatrix of $A$ for the shift of finite
type on the symbols in $\pi^{-1} \{2,3\}$ and suppose $B$
is reducible. For $i=2,3,$ let $B_{i}$ be the transition
submatrix of $A$ for the shift on the symbols in $\pi^{-1}\{i\}.$
Suppose $B_2$ and $B_3$ are the irreducible components of
$B.$ Assume that each $B_i$ is primitive or $B_i=[0]$ for each $i.$
Then there exists a saturated compensation function $\com\in C(X)$
such that $G\in C(Y)$ and, for all $\tau \in \R, \tau \com\in Bow
(X,\sigma).$
\end{theorem}
\begin{remark}\label{nyari}
\textnormal
{For the general case when the alphabet of $Y$ has more than three symbols,
we can find a measurable function $G:Y\rightarrow \R$ such that $G\circ \pi$ satisfies (\ref{cosc}) and $\sup_{n\geq 1} \var_{n}(S_{n}(G\circ\pi))<\infty.$ Such a $G$ is continuous except on a set of measure zero with respect to every $\sigma_{Y}$-invariant measure. This example could be studied in connection with Question 4 in Section 8.}\\
\end{remark}
\begin{proof}
Without loss of generality, assume that $23$ is allowable in $Y$ and $B_i$ is primitive for
$i=2,3.$ Let $B'$ be the transition matrix for the shift of finite type on symbols in $\{2,3\}$ and let $b_i$ be corresponding maximum eigenvalue.
Define $E_1$ and $E_2$ by
$$E_1=\{y\in Y: y=y_{0}\cdots y_{p-1}i^{\infty}\textnormal{ for some
} p \geq 1,y\neq i^{\infty} \textnormal{ for }i=1,2,3 \}$$ and
$E_2=\{i^{\infty}: i=1,2,3.\}.$
Define $G: Y\rightarrow \R$ by (defining $\vert\pi^{-1}[12^{0}3^{m}1]\vert=\vert\pi^{-1}[13^{m}1]\vert,$\\
\noindent $\vert\pi^{-1}[12^{0}1]\vert=\vert\pi^{-1}[13^{0}1]\vert=1) $
\begin{equation*}\label{arigatoucomfun} G(y)= \begin{cases}
 \log ({\vert \pi^{-1}[12^{n-1}3^{m}1]\vert}/{\vert \pi^{-1}[12^{n}3^{m}1]\vert
 })
& \text{
 if }y\in[2^{n}3^{m}1],n,m\geq 1\\
 \log ({\vert \pi^{-1}[1i^{n-1}1]\vert}/{\vert
\pi^{-1}[1i^{n}1]\vert}) & \text{
 if }y\in [i^{n}1], i=2,3, n\geq 1\\
\lim_{m\rightarrow \infty} \log (\vert \pi^{-1}[12^{n-1}3^m1]\vert/ {\vert
\pi^{-1}[12^{n}3^{m}1]\vert})   & \text{
 if }y =2^{n}3^{\infty}, n\geq 1 \\
 -\log b_i & \text{ if } y=i^{\infty}, i=2,3\\
0 & \text{
 if }y\in [1].
\end{cases}
\end{equation*}
Note that $G$ is continuous. We claim that for all $y \in Y \setminus E_1$ we have (\ref{upperb}) and (\ref{lowerb}).
For $y \in E_2,$ (\ref{upperb}) and (\ref{lowerb}) hold easily. Let $y
\notin E_1 \cup E_2.$ Then (\ref{lowerb}) holds by using similar arguments
as in the proof of Theorem \ref{thm1} (2)(a). To show (\ref{upperb}), notice first that by hypothesis
there are two communicating classes $C_1, C_{2}$ in
$X_B$ such that there is no
arrow in the graph of $X$ from members of $C_{1}$ to members
$C_{2}$. Note that we can write
$y=(1^{l_{1}})s_{1}1^{l_2}s_{2}1^{l_3}s_{3}\cdots,$ where each $s_i$
is a word of length $n(i)$ in $Y_{B'},$ and $l_{i}\geq 1$. Let
$s_i=y^{i}_{0}\cdots y^{i}_{n(i)-1}$ for each $i$.
Using the Perron-Frobenius
Theorem and the above fact, we find bounds for $\vert
\pi^{-1}[1s_{i}1]\vert, \vert \pi^{-1}[1y^{i}_{0}\cdots
y^{i}_{k}]\vert,$ and $\vert \pi^{-1}[1y^{i}_{k+1}\cdots
y^{i}_{n(i)-1}1]\vert$ for any $k\geq 0.$ Straightforward computations
show that $e^{(S_nG)(y)}\vert F_n(y)\vert $ is bounded uniformly.
To show that $G\circ \pi\in Bow(X,\sigma)$, we use similar arguments as in the second part of the
proof of Lemma \ref{claim1}.
\end{proof}

\section{Examples for
Theorems \ref{prop1}, \ref{prop2}, \ref{prop3}, and
 \ref{arigatou}}\label{lexamples}
We first give examples that illustrate each theorem in Section
\ref{extensions} and then present examples for which we do not know
the existence of a saturated compensation function or uniqueness of
the equilibrium state for it.
\begin{example}\label{p1e0} {An example for Theorem \ref{prop1}}
\end{example}
Let $X \subset \{1, 2, 3,4,5\}^{{\N}}$ be the shift of finite
type with the transition matrix $A$ given below. Define $\pi$ by
$\pi(1)=1, \pi(2)=\pi(3)=2,$ and $\pi(4)=\pi(5)=3.$ Let $B$ be the transition matrix of $\pi^{-1}\{2,3\}$ and
$B'$ be the
transition matrix of $Y\subset \{1, 2, 3\}^{{\N}}.$ Then
\begin{displaymath}
A= \left( \begin{array}{ccccc} 0&1&1&1&1\\
1 & 1 & 1 &1 &1  \\
1 & 0 & 1 &0 &1\\
1 & 0 & 0 & 1 & 1 \\
1 & 0 & 0 & 0 & 1
\end{array} \right)
\textnormal { and }
B'= \left( \begin{array}{ccc} 0&1&1\\
1 & 1 & 1   \\
1 & 0 & 1 \\
\end{array} \right)
\end{displaymath}

\noindent Then $h_\text{top}(X_B)=0$, $\vert \pi^{-1}[12^{n}1]\vert=\vert
\pi^{-1}[13^{n}1]\vert=n+1,$ and $\vert
\pi^{-1}[12^{k}3^{l}1]\vert=k+l+1$ for all $k,l\geq 1.$ Hypothesis I
of Theorem \ref{prop1} is satisfied. There is a compensation
function $G \circ \pi \in C(X),$ where $G\in C(Y)$ is defined by
$$G(y)= \left\{ \begin{array}{lll}
 \log ({n}/({n+1})) & \textrm{ if } y \in [2^{n}1] \cup [3^{n}1]\cup [2^{k}3^{l}1]\textnormal{ for }n\geq 1,k+l=n \\
 0 & \textrm{ if }y \in [1]\cup \{2^{\infty}, 3^{\infty}, 2^{k}3^{\infty} k\geq 1\}.
 \end{array}\right.
$$
Condition [C$'$] is satisfied, and so for each $0\leq \tau<1, \tau G\circ \pi$
has a unique equilibrium state. The unique equilibrium state $\mu$ is not Gibbs, by
the first part of the proof of Lemma \ref{claim1}. Moreover,
$(\sigma, \mu)$ is an exact endomorphism, hence strongly mixing. To
see this, use the fact that $\tau G\circ \pi$ is a grid function and apply
Theorem \ref{MPkey} to show that $\tau G\circ \pi$ satisfies the RPF
condition.

\begin{example}\label{p1e3} {An Example for Theorem \ref{prop1}}
\end{example}
This is an example for Theorem \ref{prop1} (a) but not for (b).
However, $G\circ\pi$ is still a grid function and we find a unique
equilibrium state for it.  Let $X \subset \{1, 2, 3,4,5\}^{{\N}}$
and $Y\subset \{1, 2,3\}^{{\N}}$ be the shifts of finite type
determined by the transitions given in Figure \ref{fig3}. Define
$\pi$ by $\pi(1)=1, \pi(2)=\pi(3)=2,$ and $\pi(4)=\pi(5)=3.$
\bigskip
\begin{center}
\xymatrix{ & & 2 \ar@(ul,dl) \ar@{<->}[rr]  \ar@{->}[dd]& &   1
\ar@{<->} [dd] \ar@{<->}[ddll]  \ar@{<->}[rr]& &  3 \ar@(ur,dr),
\ar@{->}[ddll]
 & &
2 \ar@(ul,dl)\ar@{<->}[rr] \ar@{->}[ddrr] & & 1 \ar@{<->}[dd]\\
& & & & & &   \ar [r]^\pi & &  \\
& & 4  \ar@(ul,dl) \ar@{->}[rr]& & 5 \ar@(dl,dr) & & & & & & 3
\ar@(ur,dr )}
\end{center}
\bigskip
\begin{figure}[h]
{} \caption{$X,Y$ and $\pi$ in Example \ref{p1e3}} \label{fig3}
\end{figure}
\noindent Then $h_\text{top}(X_B)=0$, $\vert \pi^{-1}[12^{n}1]\vert=2,
\vert \pi^{-1}[13^{n}1]\vert=n+1,$ and $\vert
\pi^{-1}[12^{k}3^{n}1]\vert=n+1$ for all $n,k.$ Hypothesis I of
Theorem \ref{prop1} is satisfied. There is a compensation function
$G \circ \pi \in C(X),$ where $G\in C(Y)$ is defined by
$$G(y)= \left\{ \begin{array}{lll}
 \log ({n}/({n+1})) & \textrm{ if } y \in[3^{n}1]\textnormal{ for }n\geq 1, \\
 0 & \textrm{ if }y \in[2^{n}]\cup [1]\cup \{3^{\infty}\}, \textnormal{ for } n\geq2,\\
 -\log 2 & \textrm{ if }y \in[21].
 \end{array}\right.
$$
Condition [C$'$] is not satisfied, so $G$ is not a grid
function associated with $\Sigma_{B'}$ and $\p$ defined as in the proof of
Theorem \ref{prop1}, but for each $0\leq \tau<1, \tau G$ is
still a grid function with a different partition. Notice that $\tau G\circ\pi$ is also a grid function.
Using Theorem \ref{MPkey}, there exists a unique equilibrium state for it. The
unique equilibrium state $\mu$ is not Gibbs, by Lemma \ref{claim1}.
By the proof of Theorem \ref{thm1}(2)(c), $(\sigma, \mu)$ is an
exact endomorphism, hence strongly mixing.
\begin{example}\label{p1e4}
\end{example}
This is an example to which we cannot apply Theorem \ref{prop1}, but
we still have the existence of a saturated compensation function and
uniqueness for it. Let $X\subset \{1, 2, 3,4,5\}^{\N}$ and $Y\subset
\{1, 2,3\}^{\N}$ be the shifts of finite type determined by the
transitions given in Figure \ref{fign}. Define the factor map $\pi$
by $\pi(1)=1, \pi(2)=\pi(3)=2,$ and $\pi(4)=\pi(5)=3.$
\bigskip
\begin{center}
\xymatrix{ & & 5 \ar@(ul,dl) \ar@{<->}[rr]  & &   1  \ar@(ul, ur)
\ar@{<->} [dd] \ar@{<->}[ddll] \ar@{<->}[ddrr] & & & &
3 \ar@(ul,dl)\ar@{<->}[rr] & & 1\ar@(ur,dr)\ar@{<->}[dd] \\
& & & & & &  \ar [r]^\pi & &  \\
& & 2  \ar@(dl,dr) \ar@{->}[rr]  \ar@{->}[uu] & & 3 \ar@(dl,dr)
\ar@{->}[rr] & &  4 \ar@(dl,dr)  & & & & 2\ar@(ur,dr)\ar@{->}[uull]
}
\end{center}
\begin{figure}[h]
{} \caption{$X,Y,$ and $\pi$ in Example \ref{p1e4}} \label{fign}
\end{figure}
\noindent Since $h_\text{top}(X_B)=0$, $\vert \pi^{-1}[12^{n}1]\vert=n+1,
\vert \pi^{-1}[13^{n}1]\vert=2,$ and $\vert
\pi^{-1}[12^{k}3^{n}1]\vert=k+1$ for all $n,k$, we have

$$\lim_{n\rightarrow\infty}\frac{\vert\pi^{-1}[12^{n-1}1]\vert}{\vert\pi^{-1}[12^{n}1]\vert}=
\lim_{n\rightarrow\infty}\frac{\vert\pi^{-1}[13^{n-1}1]\vert}{\vert\pi^{-1}[13^{n}1]\vert}=1,$$
but
$$\lim_{n\rightarrow\infty}\frac{\vert\pi^{-1}[12^{k-1}3^{n}1]\vert}{\vert\pi^{-1}[12^{k}3^{n}1]\vert}=\frac{k}{k+1}\neq
1.$$ Thus hypothesis I of Theorem \ref{prop1} is not satisfied.
Using the proof of Theorem \ref{prop1}(a), we can find $G$ that
satisfies (\ref{Nice eq}), but it is not continuous on $Y.$ Such a
$G$ is defined by
$$G(y)= \left\{ \begin{array}{lll}
 \log ({n}/({n+1})) & \textrm{ if } y \in[2^{n}1] \textnormal{ for }n\geq 1, \\
 \log ({k}/({k+1})) & \textrm{ if } y \in[2^{k}3^{n}1] \textnormal{ for }n\geq 1, k\geq 1, \\
0 & \textrm{ if }y \in [1]\cup [3]\cup \{2^{\infty},
2^{k}3^{\infty}\}, k\geq 0.
 \end{array}\right.
$$
We now modify $G$ to find a continuous saturated compensation
function.  Noting that $\vert \pi^{-1}[12^{k}3^{n}]\vert=k+1$ for
any $n,$ replacing 0 by $\log(k/(k+1))$ at $2^{k}3^{\infty},$ we get

$$\tilde G(y)= \left\{ \begin{array}{lll}
 \log ({n}/({n+1})) & \textrm{ if } y \in[2^{n}1] \cup [2^{n}3] \textnormal{ for }n\geq 1, \\
0 & \textrm{ if }y \in [1]\cup [3]\cup \{2^{\infty}\}.
 \end{array}\right.
$$
Then $\tilde G\circ \pi$ is a saturated compensation function. Since
for $0\leq \tau<1, \tau \tilde G\circ \pi$ is a grid function, using
Theorem \ref{MPkey}, there exists a unique equilibrium state for it.
The unique equilibrium state $\mu$ is not Gibbs. $(\sigma, \mu)$ is
an exact endomorphism, hence strongly mixing.

\begin{example}\label{p2e1} {An example for Theorem \ref{prop2} (and Theorem \ref{prop3}(2)) }
\end{example}
Let $X\subset \{1, 2, 3, 4,5\}^{\N}$ and $Y= \{1, 2,3\}^{\N}$ be
the shifts of finite type determined by the transitions given in
Figure \ref{fig5}, and define $\pi$ by $\pi(1)=1, \pi(2)=\pi(3)=2,$
and $\pi(4)=\pi(5)=3.$
\bigskip
\begin{center}
\xymatrix{ & & 2 \ar@(ul,dl) \ar@{<->}[rr]  & &   1  \ar@{<->} [dd]
\ar@{<->}[ddll] \ar@{<->}[ddrr] & & & &
2 \ar@(ul,dl)\ar@{<->}[rr]  & & 1 \ar@{<->}[dd] \\
& & & & & &  \ar [r]^\pi & &  \\
& & 3    \ar@{<->}[uu] & & 4 \ar@(dl,dr) \ar@{<->}[rr] & & 5  & & &
& 3 \ar@(ur,dr)}
\end{center}
\begin{figure}[h]
{} \caption{$X,Y,$ and $\pi$ in Example \ref{p2e1}} \label{fig5}
\end{figure}
\noindent Let $a=(1+\sqrt 5)/2.$ Then
$$\lim_{n\rightarrow \infty}\frac{\vert \pi^{-1}[12^{n-1}1]\vert}{\vert
\pi^{-1}[12^{n}1]\vert}=\lim_{n\rightarrow \infty}\frac{\vert
\pi^{-1}[13^{n-1}1]\vert}{\vert \pi^{-1}[13^{n}1]\vert}=\frac{1}{a},
\textnormal{ and } h_\text{top}(X_B)=\log a.$$ Define $G$ by (defining $\vert
\pi^{-1}[12^{0}1]\vert=1$)
$$G(y)= \left\{ \begin{array}{lll}
 \log ({\vert \pi^{-1}[12^{n-1}1]\vert}/{\vert
\pi^{-1}[12^{n}1]\vert}) & \textrm{ if } y \in[2^{n}1] \cup [3^{n}1]\textnormal{ for }n\geq 1, \\
 -\log a & \textrm{ if } y \in \{2^{\infty}, 3^{\infty}\},\\
0 & \textrm{ if }y \in [1].
 \end{array}\right.
$$
Since Condition [C$'$] is satisfied, $\tau G\circ \pi $ has a
unique equilibrium state for $0\leq \tau <1.$ Note that this is also
an example to illustrate Theorem \ref{prop3}(2). Thus $G\circ \pi\in Bow(X,\sigma). $

\begin{example}\label{newe}{An example for Theorem \ref{prop3}(2)}
\end{example}
Let $X\subset \{1, 2, 3, 4,5,6,7,8\}^{\N}$ be the shift of finite
type with the transition matrix $A$ given below. Define $\pi$ by
$\pi(1)=1, \pi(2)=\pi(3)=\pi(4)=\pi(5)=\pi(6)=2,$ $\pi(7)=\pi(8)=3.$
Note that $B$ is reducible and all members in two distinct
communicating classes $\{2,3\}$ and $\{4,5,6\}$ are mapped to the
same symbol.
\begin{displaymath}
A= \left( \begin{array}{cccccccc} 1&1&1&1&1&1&1&1\\
1 & 1 & 1 &0 &0 & 0 & 0 & 0 \\
1 & 1 & 1 &1 &0 & 0 & 0 & 0 \\
1 & 0 & 0 & 1 & 1 & 1 & 0 &0 \\
1 & 0 & 0 & 1 & 1& 1 & 0 & 0 \\
1 & 0 & 0 &1 &1& 1 & 0 & 0 \\
1 & 0 & 0 & 0 & 0 & 0 & 1 & 1  \\
1&0&0&0&0&0&1&0 \\

\end{array} \right).
\end{displaymath}
Then $h_\text{top}(X_{B_2})=\log 3$ and $\vert \pi^{-1}[12^{n}1]\vert=
2^{n}+3^{n} +3^{n-1}(1-(2/3)^{n-1}).$ $X_{B_3}$ is primitive, and
(1) and (2) of Theorem \ref{prop3} are satisfied for $i=2,3$. Thus
there is a saturated compensation function $G\circ\pi\in
Bow(X,\sigma)$.

\begin{example} {An example for Theorem \ref{arigatou}}
\end{example}
Let $X \subset \{1, 2, 3, 4, 5, 6 \}^{{\N}}$ be the shift of finite
type with the transition matrix $A$ given below ($A^{5}>0$).  Define
$\pi$ by $\pi(1)=1, \pi(2)=\pi(3)=2,$ and $\pi(4)=\pi(5)=\pi(6)=3.$
\begin{displaymath}
A= \left( \begin{array}{cccccc} 0 & 0 & 0 &1 &0 &1 \\
0 & 1 & 1 &0 &0& 0\\
1 & 1 & 0 &0 &0 &0 \\
0&0&0&0&1&1\\
0&1&0&0&0&1\\
0&1&0&1&0&0\\
\end{array} \right)
\end{displaymath}
There are two communicating classes for $X_B,$
$C_1=\pi^{-1}\{2\}=\{2,3\}$ and $C_2=\pi^{-1}\{3\}=\{4,5,6\}.$ For
each $i=1,2,$ the irreducible matrix corresponding to $C_i$ is
primitive. Hence there is a saturated compensation function $G\circ
\pi \in C(X), G\in C(Y),$ such that
for all $\tau \in \R, \tau G\circ\pi \in Bow(X,\sigma).$\\

Following are examples for which we do not know the existence of a
saturated compensation function or uniqueness of the equilibrium
state for $\tau G\circ \pi$ for any $\tau\neq 0.$
\begin{example}\label{p1e2}
\end{example}
Let $X\subset \{1, 2, 3, 4,5,6\}^{\N}$ and $Y= \{1, 2,3\}^{\N}$ be
the shifts of finite type determined by the transitions
given in Figure \ref{fig2}. \label{examplesforthem1}Define $\pi$ by
$\pi(1)=1, \pi(2)=\pi(3)=2,$ and $\pi(4)=\pi(5)=\pi(6)=3.$
\bigskip
\begin{center}
\xymatrix{ & & 2 \ar@(ul,dl) \ar@{<->}[rr]  & &   1  \ar@(ul, ur)
\ar@{<->} [dd] \ar@{<->}[ddll] \ar@{<->}[ddrr] \ar@{<->}[rr]& &
6\ar@(ur,dr) & &
2 \ar@(ul,dl)\ar@{<->}[rr]  & & 1\ar@(ur,dr) \ar@{<->}[dd]\\
& & & & & &   \ar [r]^\pi & &  \\
& & 3  \ar@(ul,dl)  \ar@{<-}[uu] & & 4 \ar@(dl,dr) \ar@{->}[uurr]
\ar@{->}[rr] & & 5 \ar@(dl,dr)  & & & &  3\ar@(ur,dr) }
\end{center}
\begin{figure}[h]
{} \caption{$X,Y$ and $\pi$ in Example \ref{p1e2}} \label{fig2}
\end{figure}

Since $h_\text{top}(X_B)=0$, $\vert \pi^{-1}[12^{n}1]\vert=n+1,$ and
$\vert \pi^{-1}[13^{n}1]\vert=2n+1$ for all $n,$ hypothesis I of
Theorem \ref{prop1} is satisfied, but not (b). Note also that
Theorem \ref{prop3}(1) is applicable. $G \in C(Y)$ is defined by
$$G(y)= \left\{ \begin{array}{lll}
 \log ({n}/({n+1})) & \textrm{ if } y \in[2^{n}1] \textnormal{ for }n\geq 1, \\
 \log (({2n-1})/({2n+1})) & \textrm{ if } y \in[3^{n}1] \textnormal{ for }n\geq 1, \\
0 & \textrm{ if }y \in [1]\cup \{2^{\infty}, 3^{\infty}\}.
 \end{array}\right.
$$

$G$ seems not to look like a grid function for any partition, but
$G$ is a sum of two grid functions. Thus we do not know whether for
all $\tau \in (0,1)$, $\tau G\circ\pi$ has a unique equilibrium
state or not.
\begin{example}\label{p3e4}
\end{example}
Let $X \subset \{1, 2, 3,4,5\}^{{\N}}$ and $Y\subset \{1,
2,3\}^{{\N}}$ be the shifts of finite type determined by the
transitions given in Figure \ref{fig9}, and define $\pi$ by
$\pi(1)=1, \pi(2)=\pi(3)=2,$ and $\pi(4)=\pi(5)=3.$
\bigskip
\begin{center}
\xymatrix{ & & 2 \ar@(ul,dl) \ar@{<->}[rr]  & &   1  \ar@{<->} [dd]
\ar@{<->}[ddll] \ar@{<->}[ddrr] & & & &
2 \ar@{<->}[rr] \ar@(ul,dl) & & 1 \ar@ {<->}[dd]\\
& & & & & &  \ar [r]^\pi & &  \\
& & 3  \ar@{<-}[uu] \ar@(ul,dl) & & 4 \ar@(dl,dr) \ar@{<->}[rr] & &
5 & & & & 3 \ar@(ur,dr)}
\end{center}
\begin{figure}[h]
{} \caption{$X,Y,$ and $\pi$ in Example \ref{p3e4}} \label{fig9}
\end{figure}

Let $a=({1+\sqrt 5})/2.$ Then
$$\lim_{n\rightarrow \infty}\frac{\vert
\pi^{-1}[12^{n-1}1]\vert}{\vert \pi^{-1}[12^{n}1]\vert}=1, \qquad
 \lim_{n\rightarrow \infty}\frac{\vert
\pi^{-1}[13^{n-1}1]\vert}{\vert
\pi^{-1}[13^{n}1]\vert}=\frac{1}{a},$$ and $$h_\text{top}(X_{B_2})=
0\textnormal{ and }h_\text{top}(X_{B_3})=1/a.$$ The hypothesis of (1) of
Theorem \ref{prop3} is satisfied but not the hypothesis of (2). We
define $G$ by (defining $\vert
\pi^{-1}[13^{0}1]\vert=1$)
$$G(y)= \left\{ \begin{array}{llll}
 \log ({n}/({n+1})) & \textrm{ if } y \in[2^{n}1] \textnormal{ for }n\geq 1, \\
\log(\vert \pi^{-1}[13^{n-1}1]\vert/\vert \pi^{-1}[13^{n}1]\vert)& \textrm{ if } y \in[3^{n}1] \textnormal{ for }n\geq 1, \\
 -\log a & \textrm{ if } y =\{3^\infty\}, \\
0 & \textrm{ if }y \in [1]\cup \{2^{\infty}\}.
 \end{array}\right.
$$
Then $G$ is a sum of a grid function which has a unique equilibrium
state and a function in the Bowen class. We do not know whether for $\tau \in (0,1), \tau G\circ\pi$ has a unique equilibrium state or not.
\section{Problems}
In this section, we list some questions and possible directions for
future work. An overarching question is whether the geometrical properties of a SFT-NC carpet can be
understood in a fundamentally different way, in particular without using saturated compensation functions. It would be interesting to determine completely the properties of an example for which the associated symbolic factor map does not admit a saturated compensation function.
We note that the key formula (\ref{key1}) (page
\pageref{key1}) of Shin \cite{S2}, which connects the measures maximizing the weighted entropy functional $\phi_{\alpha}$ with a saturated compensation function, was proved by using the relative variational principle with a saturated compensation function as a potential.\\

{\em Question 1.} Can an SFT-NC carpet (or an
NC-carpet corresponding to another subshift not necessarily of
finite type) have several measures of full Hausdorff
dimension? Can a measure of full Hausdorff dimension fail
to have a Bernoulli natural extension?\\

{\em Question 2.} Let $T$ be an expanding map given by a
non-diagonal matrix $A$ over $\Z$.
For a compact $T$-invariant subset of the
torus, is there an ergodic measure of full Hausdorff dimension? If
so, is it unique? What are the properties of such measure(s)?
In this case the existence of measures of full Hausdorff dimension
is not known.
Bedford \cite{B2}, and Ito and Ohtsuki \cite{IO} showed that there
is a Markov partition consisting of $\vert \det A\vert$ elements. Ito and Ohtsuki
also showed that if $\mu$ is Lebesgue measure on the torus, then the
partition is a one-sided Bernoulli partition. Also, Stolot \cite{ST} gave a
construction of ``extended Markov partitions'' for the particular
toral endomorphism given by $T(x,y)=((3x+y) \text{ mod } 1,
(x+y)\text{ mod } 1),$ and this may also be
useful.\\

{\em Question 3.}  What is the Hausdorff dimension of an SFT-NC
``sponge'' (an NC carpet on an $n$-dimensional torus with $n \geq
3$)? Is there any measure of full Hausdorff dimension? If so, is it
unique? What are the properties of the measure(s)? The Ledrappier-Young formula \cite{KP} extends to higher dimensions. If
$T$ is a toral endomorphism given by a diagonal matrix
Diag$(m_1,m_2,\cdots, m_r),$ where $m_i$ are integers with
$m_{i}<m_{i+1}$ for $i=1,2,\cdots, r-1$ and $\mu$ is an ergodic
$T$-invariant probability measure on the $r$-torus,
$$\dim_H(\mu)=\sum^{r}_{\nu=1}\frac{1}{\log{m_{\nu}}}[h(\pi_{\nu}\mu)-h(\pi_{\nu-1}\mu)],$$
where the entropy $h(\pi_{\nu}\mu)$ is with respect to the
endomorphism Diag $(m_1, \cdots, m_{\nu})$ of the $\nu$-torus, and
$h(\pi_{0}\mu)=0$ by convention \cite{KP}. The existence of measures of full
Hausdorff dimension of a compact $T$-invariant set is known \cite{KP}.

If $r=3,$ we get
\begin{equation}\label{nLY}
\dim_H (\mu)=\frac{1}{\log m_{3}}[h(\pi_{3}\mu)+\alpha
h(\pi_{2}\mu)+\beta h(\pi_{1}\mu)],
\end{equation}
where $\alpha=-1+\log_{m_2}{m_3}$ and
$\beta=-\log_{m_2}m_3+\log_{m_1}m_3$. Thus we need to find the
measure(s) that maximize the right-hand side of (\ref{nLY}). To do
it, first, using the natural Markov partition for $T,$ define $X$ to
be a symbolic representation of a compact $T$-invariant subset (an
SFT-NC sponge) of the 3-torus, $Y$ to be a symbolic representation
of the projection of the compact $T$-invariant subset of the 3-torus
to the 2-torus (the first two coordinates), and $Z$ to be a symbolic
representation of the projection of the compact $T$-invariant subset
of the 3-torus to the 1-torus (the first coordinate). Define
$\pi_{1}:(X, \sigma_{X})\rightarrow (Y, \sigma_{Y})$ to be the
projection to the first two coordinates and $\pi_{2}:(Y,
\sigma_{Y})\rightarrow (Z, \sigma_{Z})$ to be the projection to the
first coordinate. Note that $\pi_{1}$ and $\pi_{2}$ are factor maps
and so $\pi_{2}\circ \pi_{1}$ is also a factor map. Then we need to
find an ergodic shift-invariant measure on $X$ that maximizes
\begin{equation*}
h_{\mu}(\sigma_{X})+\alpha h_{\pi_{1} \mu}(\sigma_Y)+\beta
h_{\pi_{2}\pi_{1}\mu}(\sigma_{Z}).
\end{equation*}
So we have a sequence of factor maps rather than just one,
complicating the symbolic dynamics considerably.\\

{\em Question 4.} Let $X,Y$ be topological dynamical systems (for
example arbitrary subshifts) and let $\pi: X \rightarrow Y$ be a
factor map. When can we find a continuous compensation
function, i.e., a function $f: X \rightarrow \R$ satisfying
$$P_X(f+\phi\circ \pi)=P_{Y}(\phi)$$ for all $\phi\in C(Y)$?
If there is no such continuous $f,$ can we always find a measurable
one? To what extent can any such (measurable but not continuous)
function substitute for a compensation function in the above
arguments?\\

\noindent \em{Acknowlegements. }\em This paper is based on the author's Ph.D. dissertation. I would like to thank my
advisor Professor Karl Petersen for suggesting interesting problems
and for helpful advice and discussions. I am also grateful to Professor Thomas Ward for his help revising the manuscript.
Finally I would like to thank Dr. Michael Schraudner and the referee for useful suggestions and comments.

\end{document}